\documentclass[final]{siamart220329}

\usepackage{graphicx}
\usepackage{subcaption}
\usepackage{siunitx}
\usepackage{bm}
\usepackage{overpic}
\usepackage{amssymb}
\usepackage{algorithmic} 
\Crefname{ALC@unique}{Line}{Lines}

\title{A gradient flow approach for combined layout-control design of wave energy parks}
\author{Marco Gambarini\thanks{MOX, Dipartimento di Matematica, Politecnico di Milano, Piazza Leonardo da Vinci 32 20133, Milano, Italy (\email{marco.gambarini@polimi.it}, \email{gabriele.ciaramella@polimi.it}, \email{edie.miglio@polimi.it})} \and Gabriele Ciaramella\footnotemark[1] \and Edie Miglio\footnotemark[1]}

\headers{Gradient flow optimization of wave energy parks}{Marco Gambarini, Gabriele Ciaramella, and Edie Miglio}

\graphicspath{{figures/}}

\begin{document}
\maketitle

\begin{abstract}
Wave energy converters (WECs) represent an innovative technology for power generation from renewable sources (marine energy). Although there has been a great deal of research into such devices in recent decades, the power output of a single device has remained low.
Therefore, installation in parks is required for economic reasons. The optimal design problem for parks of WECs is challenging since it requires the simultaneous optimization of positions and control parameters. While the literature on this problem usually considers metaheuristic algorithms, we present a novel numerical framework based on a gradient-flow formulation. This framework is capable of solving the optimal design problem for WEC parks. In particular, we use a low-order adaptive Runge-Kutta scheme to integrate the gradient-flow equation and introduce an inexact solution procedure. Here, the tolerances of the linear solver used for projection on the constraint nullspace and of the time-advancing scheme are automatically adapted to avoid over-solving so that the method requires minimal tuning. We then provide the specific details of its application to the considered WEC problem: the goal is to maximize the average power produced by a park, subject to hydrodynamic and dynamic governing equations and to the constraints of available sea area, minimum distance between devices, and limited oscillation amplitude around the undisturbed free surface elevation. A suitable choice of the discrete models allows us to compute analytically the Jacobian of the state problem’s residual. Numerical tests with realistic parameters show that the proposed algorithm is efficient, and results of physical interest are obtained.

\end{abstract}



\section{Introduction}
This work introduces a novel numerical framework for the solution of optimal design problems for parks of wave energy converters (WEC). The problem is formulated as 
\begin{align}
\min_{\bm{w}} \quad f(\bm{w}) \quad  \text{s.t.} \quad
\bm{e}(\bm{w}) := \begin{bmatrix}
    \bm{e}_{\gamma}(\bm{w}) \\
    \bm{e}_{\zeta}(\bm{w}) 
\end{bmatrix} = 0, \quad \bm{h}(\bm{w}) := \begin{bmatrix}
    \bm{h}_{sl}(\bm{w})\\ 
    \bm{h}_{ad}(\bm{w})\\ 
    \bm{h}_{md}(\bm{w})
\end{bmatrix} \leq 0.
\label{eq:introopt}
\end{align}
Here, the cost function is $f(\bm{w}):=-P(\bm{w})$, with $P(\bm{w})$ the total power absorbed by the park. Vector $\bm{w}$ contains the decision variables, corresponding to devices' positions and control parameters, and the state variables, given by hydrodynamic coefficients $\hat{\bm{\gamma}}$ and oscillation amplitudes $\hat{\bm{\zeta}}$. Vector $\bm{e}(\bm{w})$ defines the equality constraints, corresponding to the state problem: in particular, $\bm{e}_\gamma(\bm{w})$ is the residual of the hydrodynamic equation and $\bm{e}_\zeta(\bm{w})$ is the residual of the dynamic equation. Furthermore, $\bm{h}(\bm{w})$ is the vector of inequality constraints: the slamming constraint function $\bm{h}_{sl}(\bm{w})$, requiring the devices to remain immersed as they oscillate; the admissible domain constraint function $\bm{h}_{ad}(\bm{w})$, describing the sea area available for the park; and finally, the constraint of minimum distance between each pair of devices, expressed by function $\bm{h}_{md}(\bm{w})$.

The WEC optimal design problem \cref{eq:introopt} is currently relevant for applications, because wave energy is being included in renewable energy strategies \cite{eucommdelivering}, and for it to be cost-effective, devices need to be installed in parks \cite{Goteman_2022}.
The preliminary design of a WEC park requires specifying the positions of the devices and the control strategy to be applied, in such a way that optimization criteria are satisfied. 
The final aim should be to minimize some suitable measure of the cost of energy; however, because of the difficulty of estimating costs for such an emerging technology \cite{Goeteman2020}, the average power output is often considered as the objective function instead \cite{Teixeira2022}, and this is the choice made in this work. 
Many approaches have been considered for modeling and optimizing WEC parks \cite{Golbaz2022, Yang2022}. 
The most common hydrodynamic model is based on linear potential flow theory, while optimization is typically tackled through metaheuristic algorithms. Although array layout and control strategy are usually optimized separately, it has been recognized that their simultaneous optimization can lead to significantly better results in terms of power \cite{GarciaRosa2015}. 
This approach is referred to as control co-design \cite{Ringwood2023, Neshat_2020}.

The present work introduces an adaptive gradient-flow algorithm and applies it to the simultaneous optimization of layout and control for WEC parks, with the average power as the objective function. 
To the best of our knowledge, the gradient-flow approach has not been considered for wave energy applications before, and more generally gradient-based approaches are rare in the literature on control co-design for WEC arrays.

Our optimization algorithm is based on the gradient-flow framework originally introduced in \cite{Tanabe1980}, which allows a rather natural treatment of very general constraints, provided that their Jacobian is available.
The resulting system of ordinary differential equations (ODE) is solved through an adaptive Runge-Kutta (RK) method \cite{Hairer1993}. 
To avoid over-solving while still achieving convergence within a specified threshold, an approach akin to inexact optimization algorithms \cite{Brown2021, Byrd2008, Yan1999} is adopted: the tolerances of RK stepsize adaptivity and of the conjugate gradient (CG) solver used to compute the right hand side of the ODE are dynamically determined.
For our application, we consider arrays of cylindrical heaving point absorber devices with reactive control \cite{MariaArenas2019}, i.e., the control parameters are the equivalent damping and stiffness coefficients of the electric generator, usually referred to as power take-off system (PTO). 
Hydrodynamics is modeled using linear potential flow theory \cite{Newman2018}, and irregular waves force the system. 
The wave spectrum is approximated by a superposition of a finite number of monochromatic waves.
We consider the frequency domain fluid-structure interaction model presented in \cite{Yilmaz1998, Child2010, Child2011}, based on interaction theory \cite{Kagemoto1986}. The velocity potential is written as a sum of cylindrical harmonics with unknown coefficients. 
The latter, together with the devices' oscillation amplitudes, are computed by solving a linear system. 
The system matrix can be explicitly differentiated with respect to the coordinates of the devices and to the control parameters. 
This possibility has been leveraged in previous works \cite{jacopo, Gambarini2024} to perform optimization in regular waves through an adjoint approach.


The paper is structured as follows. \Cref{sec:optalgos} is dedicated to the optimization algorithm, written for a very general class of problems. In \cref{sec:modeling}, we specialize the framework to our problem of interest, describing in detail the terms appearing in \eqref{eq:introopt} and detailing the computation of derivatives.  In \cref{sec:numtests}, several numerical experiments are presented in order to show the effectiveness of the proposed method. Finally, the limitations and some possible extensions of the present work are discussed in \cref{sec:conclusions}.

\section{Optimization algorithm}\label{sec:optalgos}
In this section, we describe the optimization algorithm used. First, in \cref{sec:gradflow} we review the gradient-flow approach. Then, in \cref{sec:scaling} we describe the scalings used to improve the conditioning of the problem. In \cref{sec:timestepping}, we finally discuss the time stepping method used and introduce an inexact solution strategy.

\subsection{Gradient flow}\label{sec:gradflow}
The gradient-flow method is a continuous version of gradient descent, where the optimal solution is sought by solving a system of ODEs whose forcing term depends on the gradient vector. An advantage of gradient-flow-like algorithms over others, such as penalty methods, is that they generate trajectories that can move along an equality constraint, instead of oscillating around it \cite{Gabriele1977}. 
For a uniform treatment of the constraints, a quadratic slack variable method \cite[Sec. 3.3.2]{Bertsekas1999} is used to recast inequality constraints as equality constraints.
To do that, vector $\bm{w}$ is augmented with slack variables $\bm{s}$, and inequality constraints are replaced by the equality constraint
\begin{equation*}
\bm{g}_\mathcal{I}(\bm{w}) := \bm{h}(\bm{w}) + \bm{s} \odot \bm{s} = 0,
\end{equation*}
where $\odot$ denotes the element-wise product.
Then, a single constraint vector $\bm{g}$ is introduced, $\bm{g}(\bm{w}) = [\bm{e}(\bm{w}), \; \bm{g}_\mathcal{I}(\bm{w})]$,
and problem \cref{eq:introopt} is restated as
\begin{equation}
\min_{\bm{w}} \quad f(\bm{w}) \quad \text{s.t.} \quad
\bm{g}(\bm{w}) = 0.
\label{eq:eqconstropt}
\end{equation}
For problem \cref{eq:eqconstropt},
we consider the constrained formulation introduced in  \cite{Tanabe1980}:
\begin{equation}
\frac{d\bm{w}}{dt} = - \underbrace{(D_w\bm{g}(\bm{w}))^+ \bm{g}(\bm{w})}_{\bm{\Psi}_g(\bm{w})} - \underbrace{\left[ I - (D_w\bm{g}(\bm{w}))^+ D_w\bm{g}(\bm{w}) \right]\nabla f(\bm{w})}_{\bm{\Psi}_f(\bm{w})}  =: \bm{\Psi}(\bm{w}).
\label{eq:gradflowtanabe}
\end{equation}
$D_w\bm{g}(\bm{w})$ is the differential of the constraint vector, 
\begin{equation*}
(D_w\bm{g}(\bm{w}))^+=(D_w\bm{g}(\bm{w}))'[D_w\bm{g}(\bm{w}) (D_w\bm{g}(\bm{w}))']^{-1}
\end{equation*} 
its pseudoinverse operator, and $(D_w\bm{g}(\bm{w}))'$ the adjoint operator of $D_w\bm{g}(\bm{w})$.
The first contribution to $\bm{\Psi}(\bm{w})$, $\bm{\Psi}_g(\bm{w})$, is orthogonal to the tangent space to the local level set of $\bm{g}(\bm{w})$, and it is sometimes referred to as range space step \cite{Feppon2020}, or feasibility step \cite{Gonzaga2004}; the second, $\bm{\Psi}_f(\bm{w})$, belongs to the tangent space and is sometimes called null space step, or optimality step.
We can rewrite $\bm{\Psi}(\bm{w})$ as
\begin{equation*}
\bm{\Psi}(\bm{w}) = -(D_w\bm{g}(\bm{w}))' [D_w\bm{g}(\bm{w}) (D_w\bm{g}(\bm{w}))']^{-1} (\bm{g}(\bm{w}) -  D_w\bm{g}(\bm{w}) \nabla f(\bm{w})) - \nabla f(\bm{w}).
\end{equation*}
Thus, the term $\bm{\Psi}(\bm{w})$ can be computed by either solving the linear system 
\begin{equation}
D_w\bm{g}(\bm{w}) (D_w\bm{g}(\bm{w}))' \bm{\Lambda} =  \bm{g}(\bm{w}) - D_w\bm{g}(\bm{w}) \nabla f(\bm{w})
\label{eq:schursys}
\end{equation}
and then substituting $\bm{\Lambda}$ in $\bm{\Psi}(\bm{w}) = -(D_w\bm{g}(\bm{w}))' \bm{\Lambda} - \nabla f(\bm{w})$, or by solving the saddle point problem
\begin{equation*}
\begin{bmatrix}
I & (D_w\bm{g}(\bm{w}))' \\
D_w\bm{g}(\bm{w}) & 0
\end{bmatrix}
\begin{bmatrix}
\bm{\Psi}(\bm{w}) \\ \bm{\Lambda}
\end{bmatrix} = 
\begin{bmatrix}
-  \nabla f(\bm{w}) \\ - \bm{g}(\bm{w})
\end{bmatrix}.
\end{equation*}
In the first case, the system matrix is symmetric positive definite, while in the second case, it is symmetric indefinite. In the following, we will consider the first approach, which allows the use of the conjugate gradient method.

For the purpose of checking convergence, we observe that $\bm{\Psi}(\bm{w})=0$ if and only if there exists a vector $\bm{\Lambda}$ such that $(D_w\bm{g}(\bm{w}))' \bm{\Lambda} + \nabla f(\bm{w}) = 0$, which is the statement of the first order optimality condition for problem \cref{eq:optprobcompact} \cite[Chap. 12]{Nocedal2006}. Therefore, $\|\bm{\Psi}(\bm{w})\|_2$ is taken as a convergence indicator.

\subsection{Scaling and preconditioning}\label{sec:scaling}
In our problem, as in most problems of practical interest, the vector $\bm{w}$ contains quantities of different scales and units. For numerical purposes, they are divided by reference values to obtain a better scaling of the problem, that is, to prevent a spurious domination of a set of variables over another \cite[Sec. 2.2]{Nocedal2006}. We define nondimensional variables $\overline{\bm{w}}$, $\overline{\bm{g}}(\overline{\bm{w}})$ and $\overline{f}(\overline{\bm{w}})$, positive and diagonal scaling matrices $S_w$ and $S_g$, and a scaling factor $f_0$ as
\begin{equation}
\bm{w} = S_w \overline{\bm{w}}, \quad 
\bm{g}(\bm{w}) = S_g \overline{\bm{g}}(\overline{\bm{w}}), \quad
f(\bm{w}) = f_0 \overline{f}(\overline{\bm{w}}).
\label{eq:scalings}
\end{equation}
Since the variable $t$ in \eqref{eq:gradflowtanabe} has no physical meaning, we consider it nondimensional. 
From the chain rule, we obtain
\begin{align}
\nabla f(\bm{w}) &= 
\frac{\partial f}{\partial \bm{w}} = f_0 S_w^{-1} \frac{\partial \overline{f}}{\partial \overline{\bm{w}}} = f_0 S_w^{-1} \overline{\nabla} \, \overline{f} (\overline{\bm{w}}) \label{eq:gradfscaled} \\
D_w\bm{g}(\bm{w}) &=  S_g D_{\overline{w}}\overline{\bm{g}}(\overline{\bm{w}})  S_w^{-1}.
\label{eq:diffgscaled}
\end{align}
The right-hand side of the ODE is then written as
\begin{equation*}
\overline{\bm{\Psi}}(\overline{\bm{w}})  = - (D_{\overline{w}}\overline{\bm{g}}(\overline{\bm{w}}) )' [D_{\overline{w}}\overline{\bm{g}}(\overline{\bm{w}}) (D_{\overline{w}}\overline{\bm{g}}(\overline{\bm{w}}) )']^{-1} (\overline{\bm{g}}(\overline{\bm{w}})  - D_{\overline{w}}\overline{\bm{g}}(\overline{\bm{w}})  \overline{\nabla}\,\overline{f}(\overline{\bm{w}}) ) -  \overline{\nabla}\,\overline{f}(\overline{\bm{w}}) .
\end{equation*}
Using \cref{eq:scalings}, \cref{eq:gradfscaled}, and \cref{eq:diffgscaled}, we get
\begin{equation*}
\begin{split}
\overline{\bm{\Psi}}(\overline{\bm{w}}) =& -S_w (D_w\bm{g}(\bm{w}))' [D_w\bm{g}(\bm{w}) S_w^2 (D_w\bm{g}(\bm{w}))']^{-1} [\bm{g}(\bm{w}) - f_0^{-1} D_w\bm{g}(\bm{w}) S_w^2 \nabla f(\bm{w})] \\&- f_0^{-1} S_w \nabla f(\bm{w}),
\end{split}
\end{equation*}
which shows that the time evolution $\overline{\bm{w}}(t)$ is affected $S_w$, but not by $S_g$. Thus, $S_g$ can then be chosen to improve the condition number of the linear system \eqref{eq:schursys} through row equilibration. Indeed, the involved linear operator is
\begin{equation*}
\Sigma = D_{\overline{w}}\overline{\bm{g}}(\overline{\bm{w}}) (D_{\overline{w}}\overline{\bm{g}}(\overline{\bm{w}}))' = S_g^{-1} [D_w\bm{g}(\bm{w}) S_w^2 (D_w\bm{g}(\bm{w}))'] S_g^{-1}.
\end{equation*}
Hence, $S_g$ can be interpreted as a symmetric diagonal preconditioner for $\Sigma$. 
We choose to perform row equilibration by normalizing the rows of the representative matrix $J_g$ of $D_w\bm{g}$ with respect to the 2-norm.
If the rows of $J_g$ formed an orthogonal set, then our choice would yield $\Sigma=I$. In general, it is only guaranteed that the diagonal terms of $J_g J_g^T$ will be all ones. 
It can be proved \cite[Thm. 4.1]{Sluis1969} that the spectral condition number obtained from row equilibration $\kappa_2^{re}$ and the minimum one obtainable with a symmetric diagonal preconditioner $\kappa_2^{min}$ satisfy $\kappa_2^{re}\leq n\kappa_2^{min}$, where $n$ is the dimension of all involved matrices.
While this result indicates that there may be room for improvement, since the condition number may increase proportionally to $n$, we have observed that the proposed preconditioner yields reasonable performance in terms of number of iterations and computational time for our purpose.

\subsection{Adaptive time stepping}
\label{sec:timestepping}
We now deal with the time-discretization of \cref{eq:gradflowtanabe}. Since this is a system of ODEs, but it also originates from an optimization problem, a suitable combination of numerical approaches from the two fields can be adopted for an efficient and accurate solution.
On the one hand, there is the possibility of using numerical methods for the discretization of ODEs, possibly adapted to the optimization nature of the original problem. Such a choice was made in \cite{yamashita1980differential}, where a fourth-order Runge-Kutta scheme with an ad-hoc stepsize choice based on constraint violation was employed. On the other hand, one could rely on constrained optimization concepts, such as merit functions and filters \cite[sect. 15.4]{Nocedal2006}. This choice was made in \cite{Feppon2020}, where a merit function was used.
Our approach is similar to the former: we start from an adaptive Runge-Kutta method, and introduce criteria for the dynamic update of solver tolerances.
In all cases, to avoid the need of subiterating (and the possible need of computing higher order derivatives, if a Newton algorithm is used for subiterations), explicit algorithms are preferable.

When standard methods for ODEs are chosen, a compromise needs to be sought between number of required function evaluations, accuracy and stability. 
To limit the discussion, we consider single step, multistage methods. The number of stages can be increased with the aim of either improving the accuracy or enlarging the stability region and thus allowing larger steps. 
Adaptive time stepping can be performed by combining two methods of different orders \cite[sect. II.4]{Hairer1993}. 
Starting from iterate $\bm{w}^j$, the solution at the next step is computed using the low-order method, obtaining $\bm{w}^{j+1}$, and the high-order method, obtaining $\widehat{\bm{w}}^{j+1}$. Then, the normalized error between the two is computed as \cite[eq. 4.11]{Hairer1993}
\begin{equation}
\epsilon = \sqrt{\frac{1}{N_w} \sum_{\ell=1}^{N_w} \left(\frac{w_\ell^{j+1} - \widehat{w}_\ell^{j+1}}{\varsigma_\ell}\right)^2 }, 
\label{eq:erradaptode}
\end{equation}
where $\varsigma_\ell = \tau_{RK}^{abs} + \text{max}(|w_\ell^{j}, w_\ell^{j+1}|)\tau_{RK}^{rel}$, $\tau_{RK}^{abs}$ being the absolute tolerance and $\tau_{RK}^{rel}$ the relative tolerance. The estimated time step is then updated as
\begin{equation*}
\Delta t^{j+1}_\text{est} = (1/\epsilon)^{1/r} \Delta t^j,
\end{equation*}
where $r$ is the order of the local truncation error of the low-order method (2 for Euler's method). In practice, the variation of the time step between successive iterates is limited by defining safety factors; this leads to the actual updated time step $\Delta t^{j+1}$.
The cheapest adaptive method built in such a way in terms of stages is the combination of the 1-stage, first-order explicit Euler method and the 2-stage, second-order Heun method. If the solution is updated using the Euler method, a single evaluation of $\bm{\Psi}(\bm{w})$ is required at each time step, since the second stage is recycled as the first stage of the next step. The result is \cref{alg:eeheun}. 
\begin{algorithm}[h!]
\caption{First order Euler-Heun method}
\label{alg:eeheun}
\begin{algorithmic}[1]
\STATE Set $\bm{w}^0$, $\Delta t^0$, $\tau_\Psi$, $t_{max}$
\STATE Set $t=0$, $j=0$
\STATE $k_2 = \bm{\Psi}(\bm{w}^0)$
\WHILE{{$\|\bm{\Psi}(\bm{w}^j)\|_2 > \tau_\Psi$} \AND {$t < t_{max}$}}
	\STATE $k_1 = k_2$
	\STATE $k_2 = \bm{\Psi}(\bm{w}^j + k_1 \Delta t^j)$
	\STATE $\bm{w}^{j+1} = \bm{w}^j + k_1 \Delta t^j$
	\STATE $\widehat{\bm{w}}^{j+1} = \bm{w}^j + \frac{1}{2} (k_1 + k_2) \Delta t^j$
	\STATE Compute the relative error between the two updates using 	\eqref{eq:erradaptode} 	
	\STATE Compute the updated time step $\Delta t^{j+1}$
	\IF{$\epsilon \leq 1$}
		\STATE Stepsize accepted: $t \leftarrow t + \Delta t^j$
		\STATE $j \leftarrow j+1$
	\ENDIF
\ENDWHILE
\end{algorithmic}
\end{algorithm}

In the following, we present criteria to select the tolerance of the Runge-Kutta method $\tau_{RK}$ and the absolute tolerance on the residual of the CG solver $\tau_{CG}$ automatically and adaptively to achieve convergence within a specified tolerance on the final value of $\|\bm{\Psi}(\bm{w})\|_2$, $\tau_\Psi$. The first criterion is based on stability, while the second is based on accuracy.

The condition $\Delta t < 2$ was indicated in \cite{Feppon2020} as a stability requirement, with the following explanation. Consider the expansion
\begin{equation*}
\bm{g}(\bm{w}^{j+1}) = \bm{g}(\bm{w}^j + \Delta t \bm{\Psi}(\bm{w}^j)) = \bm{g}(\bm{w}^j) + \Delta t J_g(\bm{w}^j) \bm{\Psi}(\bm{u}^j) + \bm{R}_2,
\end{equation*}
where $\bm{R}_2$ is the remainder of the series, $\|\bm{R}_2\|=\mathcal{O}(\Delta t^2)$. 
From $J_g \bm{\Psi} = -\bm{g}$, which is a consequence of \cref{eq:gradflowtanabe}, and the triangle inequality, we obtain
\begin{equation*}
\|\bm{g}(\bm{w}^{j+1})\| \leq |1 - \Delta t| \|\bm{g}(\bm{w}^j)\| + \|\bm{R}_2\|,
\end{equation*}
and the condition stems from requiring $|1 - \Delta t|< 1$, which corresponds to a reduction of the constraint vector norm.
However, if system \eqref{eq:schursys} is solved with finite accuracy, then the maximum time step satisfying the stability requirement is reduced, as seen in the following.
Let us call $\widetilde{\Lambda}$ the numerical solution of \eqref{eq:schursys} obtained by using an iterative method, and let $\widetilde{\bm{\Psi}}(\bm{w}^j)$ be the resulting approximation of $\bm{\Psi}(\bm{w}^j)$. Then,
\begin{equation}
\bm{g}(\bm{w}^j + \Delta t \widetilde{\bm{\Psi}}(\bm{w}^j)) = \bm{g}(\bm{w}^j) + \Delta t J_g(\bm{w}^j) \widetilde{\bm{\Psi}}(\bm{w}^j) + \bm{R}_2.
\label{eq:gnextstep}
\end{equation}
From \eqref{eq:schursys}, and defining $\bm{e}_{\Lambda} = \widetilde{\bm{\Lambda}} - \bm{\Lambda}$, we have
\begin{equation}
\widetilde{\bm{\Psi}}(\bm{w}^j) = - \nabla f(\bm{w}^j) + J_g^T(\bm{w}^j)(\bm{\Lambda} + \bm{e}_\Lambda) = \bm{\Psi}(\bm{w}^j) + J_g^T(\bm{w}^j)\bm{e}_\Lambda.
\label{eq:psitildelinsol}
\end{equation}
Now, substituting \cref{eq:psitildelinsol} into \cref{eq:gnextstep} and using again $J_g \bm{\Psi} = -\bm{g}$, we get
\begin{equation*}
    \bm{g}(\bm{w}^j + \Delta t \widetilde{\bm{\Psi}}(\bm{w}^j)) = (1-\Delta t)\bm{g}(\bm{w}^j) + \Delta t J_g(\bm{w}^j)J_g^T(\bm{w}^j) \bm{e}_\Lambda + \bm{R}_2.
\end{equation*}
The last expression can be rewritten in terms of the residual of system \cref{eq:schursys} by recalling that, given any system $A\bm{x}=\bm{b}$, the relationship between the error $\bm{e}$ and the residual $\bm{r}$ is $A\bm{e}=\bm{r}$. In our case, we have $J_g(\bm{w}^j)J_g^T(\bm{w}^j) \bm{e}_\Lambda = \bm{r}$, which leads to
\begin{equation*}
\bm{g}(\bm{w}^j + \Delta t \widetilde{\bm{\Psi}}(\bm{w}^j)) = (1-\Delta t)\bm{g}(\bm{w}^j) + \Delta t \bm{r} + \bm{R}_2.
\end{equation*} 
If $\Delta t<2$, using the triangle inequality we obtain
\begin{equation*}
\|\bm{g}(\bm{w}^j + \Delta t \widetilde{\bm{\Psi}}(\bm{w}^j))\| \leq (1-\Delta t) \|\bm{g}(\bm{w}^j)\| + \Delta t\|\bm{r}\| + \|\bm{R}_2\|.
\end{equation*}
Rearranging yields the more readily interpretable inequality
\begin{equation}
\|\bm{g}(\bm{w}^j)\| - \|\bm{g}(\bm{w}^j + \Delta t \widetilde{\bm{\Psi}}(\bm{w}^j))\| \geq \Delta t\left( \|\bm{g}(\bm{w}^j)\| - \|\bm{r}\| \right) - \|\bm{R}_2\|.
\label{eq:gdescent1step}
\end{equation}
As long as $\|\bm{r} \| < \|\bm{g}(\bm{w}^j)\|$, the first-order term dominates on the second-order term and a reduction of the norm of $\bm{g}$ is guaranteed for a sufficiently small time step. However, if $\|\bm{r}\|$ approaches or exceeds $\|\bm{g}(\bm{w}^j)\|$, the first order term becomes negligible and eventually negative, and a nonincreasing value of the constraint norm may not be guaranteed for any value of $\Delta t$. We thus need to set $\|\bm{r}\| < \tau_{CG} \ll \|\bm{g}(\bm{w}^j)\|$.

Let us now focus on accuracy and assume that the exact gradient-flow trajectory $\bm{w}(t)$ converges to a constrained optimum, i.e. for any value of tolerance $\tau >0$ there exists $T>0$ such that $\|\bm{\Psi}(\bm{w}(t))\|<\tau, \; \forall t>T.$ We study the effect of the limited accuracy of the ODE and linear system solvers on the possibility of converging to within the prescribed tolerance. We denote as $\hat{\bm{w}}$ the numerical solution of ODE \eqref{eq:gradflowtanabe}. At a certain step $\hat{\bm{w}}$ along the numerically integrated trajectory, the numerical approximation of the flow map can be split into the following contributions:
\begin{equation*}
\widetilde{\bm{\Psi}}(\hat{\bm{w}}) = \bm{\Psi}(\bm{w}) + \underbrace{\left(\bm{\Psi}(\hat{\bm{w}}) - \bm{\Psi}(\bm{w})\right)}_{\text{ODE solver error}} + \underbrace{\left(\widetilde{\bm{\Psi}}(\hat{\bm{w}}) - \bm{\Psi}(\hat{\bm{w}})\right)}_{\text{linear solver error}}.
\end{equation*}
The triangle inequality leads to
\begin{equation*}
\|\widetilde{\bm{\Psi}}(\hat{\bm{w}})\| \leq \|\bm{\Psi}(\bm{w})\| + \|\bm{\Psi}(\hat{\bm{w}}) - \bm{\Psi}(\bm{w})\| + \|\widetilde{\bm{\Psi}}(\hat{\bm{w}}) - \bm{\Psi}(\hat{\bm{w}})\|.
\end{equation*}
In order to reach and observe convergence up to a tolerance, the two error terms need to be controlled. One may for instance aim to ensure that they are at least an order of magnitude smaller than $\|\bm{\Psi}(\bm{w})\|$. For the ODE solver error, we assume Lipschitz continuity of $\bm{\Psi}$ and write
\begin{equation*}
\|\bm{\Psi}(\hat{\bm{w}}) - \bm{\Psi}(\bm{w})\| \leq L \|\hat{\bm{w}} - \bm{w}\|.
\end{equation*}
For the linear solver error, we have
\begin{equation*}
\|\widetilde{\bm{\Psi}}(\hat{\bm{w}}) - \bm{\Psi}(\hat{\bm{w}})\| \leq \|J_g^+(\hat{\bm{w}})\| \|\bm{r}\|.
\end{equation*}
 The problem is now to find order of magnitude estimates of $L$ and $\|J_g^+(\hat{\bm{w}})\|$ that are cheaply computable. We adopt the following:
\begin{equation*}
L \sim \frac{\|\bm{\Psi}(\bm{w}^j) - \bm{\Psi}(\bm{w}^{j-1})\|}{\|\bm{w}^j - \bm{w}^{j-1}\|}, \quad \|J_g^+(\hat{\bm{w}})\| \sim \frac{\|J_g^T(\hat{\bm{w}}) \bm{\Lambda}\|}{\| \bm{g}(\hat{\bm{w}}) - J_g(\hat{\bm{w}}) \nabla f(\hat{\bm{w}})\|} .
\end{equation*}
All the vectors involved are intermediate quantities needed for the solution of the problem; the only additional work required to build the estimates is the computation of the norms. The tolerance of the CG solver can then be set as
\begin{equation*}
\tau_{CG}^j = k_\tau \min\left( \frac{\|\bm{\Psi}(\bm{w}^{j-1})\|\| \bm{g}(\bm{w}^{j-1}) - J_g(\bm{w}^{j-1}) \nabla f(\bm{w}^{j-1})\|}{\|J_g^T(\bm{w}^{j-1}) \bm{\Lambda}\|}, \|\bm{g}(\bm{w}^{j})\| \right),
\end{equation*}
where the stability condition implied by \eqref{eq:gdescent1step} is taken into account. The tolerance of the ODE solver can be set as
\begin{equation*}
\tau_{RK}^j = \min \left( \tau_{RK}^{j-1}, k_\tau \frac{\|\bm{\Psi}(\bm{w}^j)\| \|\bm{w}^j - \bm{w}^{j-1}\|}{\| \bm{\Psi}(\bm{w}^j) - \bm{\Psi}(\bm{w}^{j-1}) \|} \right),
\end{equation*}
where, for consistency, the error estimate used for stepsize selection needs to be computed in the euclidean norm, that is, \eqref{eq:erradaptode} is replaced with 
\begin{equation*}
\epsilon = \|\bm{w}^{j+1} - \hat{\bm{w}}^{j+1}\|/\tau_{RK}^j.
\end{equation*}
Moreover, the previous tolerance $\tau_{RK}^{j-1}$ is used as an upper bound for $\tau_{RK}^{j}$ in order to avoid large oscillations of the stepsize and frequent stepsize rejections, which may reduce the efficiency of the method. One may set $k_\tau =0.1$.
As an additional measure to prevent over-solving, the time step is halved if conjugate gradient does not converge to the required tolerance within a prescribed maximum number of iterations.

In the following, the presented algorithm is applied to the problem of optimizing positions and control parameters for wave energy converter parks.

\section{Modeling of the WEC optimal design problem}\label{sec:modeling}
In this section, we describe in detail the terms of 
\cref{eq:introopt}. Cost function and decision variables are defined in \cref{sec:costfun}.
The fluid-structure interaction model, whose residual defines the equality constraint, is presented in \cref{sec:stateprob}. The treatment of the inequality constraints, namely the slamming, available domain, and minimum distance constraints, is described in \cref{sec:ineqconstr}. In \cref{sec:derivatives}, we report the computation of the required differentials.

\subsection{Cost function and decision variables}\label{sec:costfun}
We aim to maximize the time-averaged power of an array of $N_b$ WEC devices.
The decision variables to be determined through optimization are the projections on the horizontal plane of the devices' centers, collected in vector $\bm{X} = [\bm{x}_1, \dots, \bm{x}_{N_b}] \in \mathbb{R}^{2N_b}$, and the coefficients of the PTO mechanisms: dampings $\bm{c} = [c_1, \dots, c_{N_b}]$ and stiffnesses $\bm{\kappa} = [\kappa_1, \dots, \kappa_{N_b}]$. For convenience, we define the vector of decision variables $\bm{u} = [\bm{X}, \bm{c}, \bm{\kappa}]$. We instead consider fixed the number of devices $N_b$ and their geometries, assumed to be all identical and defined by radius $R$ and draft $d$. 
We also fix the sea state, defined by significant wave height $H_s$, energy period $T_e$ and wave direction $\beta$, together with the sea depth $D$, assumed uniform.
The power is written as
\begin{equation*}
P = \lim_{T \to\infty} \frac{1}{T} \sum_{\ell=1}^{N_b} \int_0^T c_\ell \dot{\zeta}_\ell^2(t) \, dt,
\end{equation*}
$\zeta_\ell(t)$ being the elevation of the barycenter of the $\ell$-th device with respect to its hydrostatic equilibrium position.

Since the fluid-structure interaction problem is based on a linear model, detailed in \cref{sec:stateprob}, it can be formulated in the frequency domain. 
Thus, $\zeta_\ell(t)$ can be written as a linear combination of a set of harmonics. In principle, an infinite number of frequency components should be taken into account; for numerical computations, we instead consider a finite number $N_f$:
\begin{equation}
\zeta_\ell(t) = \sum_{q=1}^{N_f} \Re\left[\hat{\zeta}_{\ell q} \exp[-i(\omega_q t + \varphi_q)] \right].
\label{eq:zdecomp}
\end{equation}
Here, $\hat{\zeta}_{\ell q} \in \mathbb{C}$ is the $q$-th frequency component of the oscillation amplitude of the $\ell$-th device. The $N_f$ angular frequencies are determined through a suitable discretization of the energy spectral density of the wave climate, to be described in the next section.
In \cref{eq:zdecomp}, $\varphi_q$ is the phase of the $q$-th wave component, which does not influence the average power, and will not be considered in the following.

Thanks to the orthogonality of harmonics, the average power can be approximated as a sum of $N_f$ contributions.
\begin{equation}
P \approx \frac{1}{2} \sum_{\ell=1}^{N_b} c_\ell \sum_{q=1}^{N_f} \left( \omega_q |\hat{\zeta}_{\ell q}|\right)^2.
\label{eq:costfunction}
\end{equation}

\subsection{State problem: dynamic and hydrodynamic equations}\label{sec:stateprob}
The state problem is a linearized fluid-structure interaction problem. Thanks to linearity, it can be solved in the frequency domain, with each frequency component evolving independently from all others. First, the determination of the frequency components to be used is discussed, based on a given sea state; then, the problem to be solved for each component is introduced.

Sea states are typically described through frequency-directional spectrum functions. For simplicity, we consider unidirectional waves, and thus only concern ourselves with the  frequency dependence. 
The frequency distribution of the energy density per unit area of sea surface $E$ due to wave motion is described by the energy spectrum function $S(\omega)$ \cite{Falnes2020}: $E = \rho g \int_0^\infty S(\omega) d\omega$,
where $\rho$ is the water density and $g$ is the gravity field. Analytical forms of $S(\omega)$, parameterized in significant wave height $H_s$ and energy period $T_e$, are available in the literature. In this work, we consider the Pierson-Moskowitz spectrum \cite{Falnes2020}.
The aim of the discretization is to approximate the free-surface elevation field $\eta(x, y, t)$ as a sum of regular waves of heights $H_q$ and angular frequencies $\omega_q$, $q=1,\dots,N_f$, corresponding to wavenumbers $k_q$. 
In particular, the undisturbed free surface elevation at the center of the $\ell$-th device is 
\begin{equation}
\eta_\ell(t) = \sum_{q=1}^{N_f} \Re\left[ \hat{\eta}_{\ell q} \exp(- i\omega_q t + \varphi_q) \right], \text{ with } \hat{\eta}_{\ell q} = i \frac{H_q}{2} \exp[i k_q (x_\ell \cos\beta + y_\ell \sin \beta)].
\label{eq:etadecomp}
\end{equation}
 The energy density of the $q$-th harmonic wave component is $E_q = \rho g H_q^2/8$.   
A review of the possible strategies for discretizing wave spectra as superpositions of harmonics is reported in \cite{Chakrabarti1987}.
In this work, we proceed as follows. 
A bounded interval of frequency $(\omega_L, \omega_R)$, chosen neglecting a given fraction  of the total energy, is split into subintervals of equal width.
For each subinterval $I_q$, a wave component is defined, with angular frequency corresponding to the center of the subinterval and height $H_q$ such that $E_q = \rho g \int_{I_q} S(\omega) d\omega$.


For the fluid-structure interaction problem, we use linear potential flow theory \cite{Newman2018}.
Each device behaves as a harmonic oscillator forced by wave loads: this is described in the frequency domain by
\begin{equation}
\left( -\omega_q^2 m - i \omega_q c_\ell + \kappa_h + \kappa_\ell \right) \hat{\zeta}_{\ell q} = \hat{F}_{\ell q}.
\label{eq:dynsimple}
\end{equation}
Here, $m$ and $\kappa_h$ are mass and hydrostatic stiffness, respectively, and they are the same for all devices. 
$\hat{F}_{\ell q}$ is the $q$-th frequency component of the vertical hydrodynamic force on device $\ell$, computed from the corresponding velocity potential $\hat{\phi}_q$ as
\begin{equation}
\hat{F}_{\ell q} = i \omega \rho \int_{\Gamma_{b\ell}} \hat{\phi}_q d\Gamma,
\label{eq:forcedef}
\end{equation}
where $\Gamma_{b\ell}$ is the bottom surface of device $\ell$. 
In turn, the potential $\hat{\phi}_q$ satisfies
\begin{equation}
\begin{cases}
\Delta \hat\phi_q = 0 & \text{in $\Omega$}\\[.5em]
\dfrac{\partial \hat\phi_q}{\partial n}  =0 & \text{on $\Gamma_b$}\\[.5em]
\dfrac{\partial \hat\phi_q}{\partial z} - \dfrac{\omega^2}{g} \hat\phi_q = 0 & \text{on $\Gamma_s$}\\[.5em]
\dfrac{\partial \hat\phi}{\partial n} = -\dfrac{\partial \hat\phi_q^0}{\partial n} - i \omega_q \hat{\zeta}_{\ell q} n_z & \text{on $\Gamma_{d,\ell}$, $\ell = 1, \dots, N_b$} ,
\end{cases}
\label{eq:stateprobcontinuous}
\end{equation}
where $\Omega$ is a domain bounded by the sea bottom $\Gamma_b$, the undisturbed free surface $\Gamma_s$, the devices' surfaces $\Gamma_{d,\ell}$, and horizontally unbounded. $\hat\phi_q^0$ is the potential of an undisturbed ambient wave field corresponding to the $q$-th component of the discretization of the sea state.

In the rest of this section, we present a numerical model, derived in \cite{Child2010} and valid for cylindrical objects, for the solution of \cref{eq:stateprobcontinuous}. 
Cylindrical symmetry is leveraged to write the potential as a truncated series of cylindrical harmonics in the local coordinate system of each device. 
The $q$-th frequency component of the potential scattered by the $\ell$-th body is written as a linear combination of suitable basis functions, with unknown coefficients $\hat{\bm{\gamma}}_\ell$. 
The number of coefficients depends on parameters $N_n$ and $N_m$, which are the numbers of considered progressive and evanescent wave modes, respectively. We have $\hat{\bm{\gamma}}_\ell \in \mathbb{C}^{(2N_n + 1)(N_m+1)}$ \cite{Child2011}. These coefficients, together with the devices' oscillation amplitudes, are the unknowns of the problem, which can be written as the system
\begin{equation}
\begin{bmatrix}
M_{\gamma\gamma, q} (\bm{u}) & M_{\gamma\zeta, q} (\bm{u}) \\
M_{\zeta\gamma, q} (\bm{u}) & M_{\zeta\zeta, q}(\bm{u}) 
\end{bmatrix}
\begin{bmatrix}
\hat{\bm{\gamma}}_q \\
\hat{\bm{\zeta}}_q
\end{bmatrix} = 
\begin{bmatrix}
\bm{h}_{\gamma, q}(\bm{u}) \\
\bm{h}_{\zeta, q}(\bm{u})
\end{bmatrix}.
\label{eq:stateeq}
\end{equation}
The residuals of all such systems, for $q=1, \dots, N_f$, form the equality constraint vector of our optimization problem:
\begin{equation}
\bm{e}_q(\bm{u}, \hat{\bm{\gamma}_q}, \hat{\bm{\zeta}}_q) = \begin{bmatrix}
\bm{e}_{\gamma, q} \\
\bm{e}_{\zeta, q}
\end{bmatrix} = 
\begin{bmatrix}
M_{\gamma \gamma, q}(\bm{u}) \hat{\bm{\gamma}}_q + M_{\gamma \zeta, q}(\bm{u}) \hat{\bm{\zeta}}_q - \bm{h}_{\gamma, q}(\bm{u})\\ 
M_{\zeta \gamma, q}(\bm{u}) \hat{\bm{\gamma}}_q + M_{\zeta \zeta, q}(\bm{u}) \hat{\bm{\zeta}}_q - \bm{h}_{\zeta, q}(\bm{u})
\end{bmatrix}. \label{eq:stateresidual}
\end{equation}
In the following, we will provide the expressions of the system's blocks. We consider a single frequency and drop subscript $q$ for ease of notation.

The first block row of system \cref{eq:stateeq} corresponds to the hydrodynamic equations, which are derived from the conditions of impermeability on the surface of each body. Their expressions are
\begin{equation}
\sum_{\substack{m=1\\m\neq \ell}}^{N_b} \underbrace{B_\ell T_{m\ell}^T}_{M_{\gamma\gamma, \ell m}} \hat{\bm{\gamma}}_m - \hat{\bm{\gamma}}_\ell + \sum_{\substack{m=1\\m\neq \ell}}^{N_b} \underbrace{B_\ell T_{m\ell}^T \bm{R}_m}_{M_{\gamma\zeta, \ell m}} \hat{\zeta}_m = \underbrace{- \frac{H}{2}B_\ell \bm{a}_\ell}_{\bm{h}_{\gamma,\ell}}, \quad \ell=1,\dots,N_b,
\label{eq:hydoper}
\end{equation}
where vector $\bm{R}_m$ and matrix $B_\ell$ depend on the device's geometry but not on its position, while matrix $T_{m\ell}$ depends on the positions of bodies $m$ and $\ell$. It is the coordinate transformation matrix. $H$ is the height of the considered incident wave component, and $\bm{a}_\ell$ is a vector of ambient wave coefficients. The left-hand side of \cref{eq:hydoper} describes hydrodynamic interactions between bodies, while the right-hand side describes the external forcing due to ambient waves. 

The vertical motion of each cylinder is determined by the wave force acting on its bottom surface, given by \cref{eq:forcedef}. In the region below the $\ell$-th body, the potential is written as a linear combination of basis functions $\tilde{\bm{\Psi}}_\ell^D$ as follows:
\begin{equation}
\begin{split}
\hat{\phi}_\ell(\bm{x}_\ell^C) = \frac{g}{\omega} \Bigl\{ \Bigl[ &\bm{a}_\ell^T + \sum_{\substack{m=1\\m\neq \ell}}^{N_b} (\hat{\bm{\gamma}}_m + \hat{\zeta}_m \bm{R}_m)^T T_{m \ell} \Bigr] \tilde{B}_\ell^T \tilde{\bm{\Psi}}_\ell^D(\bm{x}_\ell^C)  \\ &  + \hat{\zeta}_\ell \Bigl(\tilde{R}_\ell^p(\bm{x}_\ell^C) + \tilde{\bm{R}}_\ell^T \tilde{\bm{\Psi}}_\ell^D(\bm{x}_\ell^C) \Bigr)  \Bigr\},
\end{split}
\label{eq:potentialbelow}
\end{equation}
where $\bm{x}_\ell^C$ is the set of local cylindrical coordinates.
Here, matrix $\tilde{B}_\ell$ and scalar function $\tilde{R}_\ell^p$ encode hydrodynamic properties of a single device, depending on its geometry, but not on its position. 
Substitution of \cref{eq:potentialbelow} into the expression of $\hat{F}_{\ell q}$ \cref{eq:forcedef} results in the following form of the dynamic equation \cref{eq:dynsimple}:
\begin{equation}
\sum_{\substack{m=1\\m\neq \ell}}^{N_b} \biggl[ \underbrace{ \frac{1}{W_\ell} (\widetilde{\bm{Y}}_\ell^D)^T \widetilde{B}_\ell T_{m\ell}^T   }_{M_{\zeta \gamma, \ell m}} \hat{\bm{\gamma}}_m +  \underbrace{\frac{1}{W_\ell} \bm{R}_m^T T_{m\ell} \widetilde{B}_\ell^T \widetilde{\bm{Y}}_\ell^D}_{M_{\zeta\zeta, \ell m}}  \hat{\zeta}_m \biggr] + \hat{\zeta}_\ell = \underbrace{- \frac{H}{2W_\ell} \bm{a}_\ell^T \widetilde{B}_\ell^T \widetilde{\bm{Y}}_\ell^D}_{\bm{h}_{\zeta, \ell}},
\label{eq:dynoper}
\end{equation}
$\ell=1,\dots,N_b$, giving the second block row of \cref{eq:stateeq}.
Here, $W_\ell$  is the mechanical impedance
\begin{equation}
W_\ell = \tilde{Y}_\ell^R - i (\omega^2m_\ell  + i\omega c_\ell - \kappa_{h,\ell} - \kappa_\ell)/(\rho g),
\label{eq:mechimped}
\end{equation}
$\tilde{\bm{Y}}_\ell^D$ is the vector of integrals of functions $\tilde{\bm{\Psi}}_\ell^D$, and $\tilde{Y}_\ell^R$ is a radiation quantity containing the integral of function $\tilde{R}_\ell^p$. 

For our purposes, the most important features of this formulation are the following. 
Firstly, the residual of the state problem can be differentiated analytically with respect to the decision variables, as detailed in \cref{sec:derivatives}. 
In particular, the dependence on the devices' coordinates is only contained in coordinate transformation matrices $T_{m\ell}$ and ambient wave vectors $\bm{a}_\ell$, the dependence on stiffness and damping is only contained in the mechanical impedances $W_\ell$, and both have differentiable analytic expressions. 
Secondly, all other quantities are independent of the decision variables and can be precomputed, so that, as vector $\bm{u}$ is updated at each optimization iteration, the residual of the state problem and its differential can be computed efficiently. 
This is in contrast with other popular modelling approaches such as the boundary element method, in which grids and interaction matrices would need to be completely reassembled as the devices' positions change.

\subsection{Inequality constraints}\label{sec:ineqconstr}
In this section, we first detail the slamming, available domain and minimal distance inequality constraints, and then convert them into equality constraints through a slack variable formulation, as explained for a generic optimization problem in \cref{sec:gradflow}.

\subsubsection{Slamming constraint}\label{sec:slamconstr}
The model presented in the previous section is linear; hence, it poses no limits on the magnitude of the oscillation amplitudes of devices. 
To ensure that the latter remain physically meaningful, in the sense that the devices do not leave the water, the slamming constraint can be enforced \cite{Backer2010}. 
In the time domain, one could impose the constraint as $\zeta_\ell(t) - \eta_\ell(t) < d$ $\forall t$, $\ell=1, \dots, N_b$.
Since we are working in the frequency domain, the root mean square value of $\zeta_\ell(t) - \eta_\ell(t)$ can instead be limited:
\begin{equation}
h_{sl, \ell} = \sum_{q=1}^{N_f} |\hat{\zeta}_{\ell q} - \hat{\eta}_{\ell q}|^2 - 2 \alpha^2 d^2 \leq 0, \quad \ell=1,\dots,N_b, \label{eq:slamconstr}
\end{equation}
where $\hat{\eta}_{\ell q}$ is the complex amplitude of the $q$-th wave component evaluated at the center of the $\ell$-th device. 
This formulation does not guarantee that the limit is never exceeded, but it rather has a statistical meaning; in particular, $\alpha$ is the parameter controlling the admissible exceedance probability \cite{Gambarini2023}.

\subsubsection{Layout constraints}\label{sec:domconstr}
In this section, the treatment of the constraints on the devices' layout is described: they are the available domain constraint and the minimum distance constraint. 

\begin{figure}[h!]
\centering
\begin{overpic}[width=0.3\textwidth]{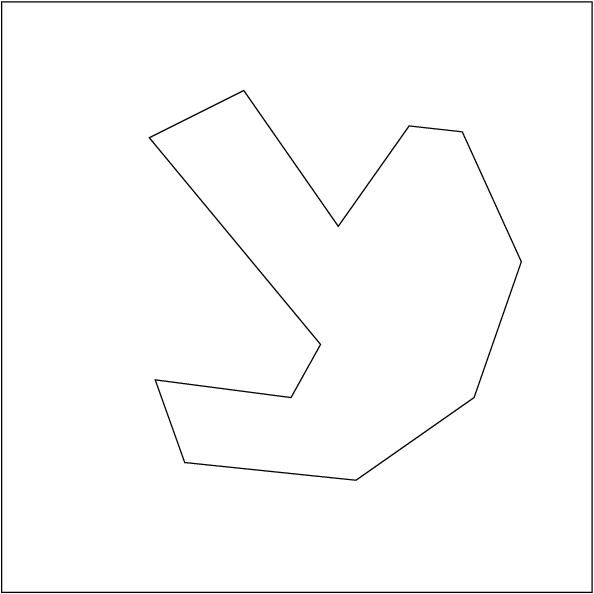}
\put (64, 48) {$\Omega_{ad}$}
\put (14, 48) {$\Omega_{ext}$}
\put (46, 78) {$\Gamma_{ad}$}
\put (-16, 20) {$\Gamma_{ext}$}
\end{overpic}
\caption{Poisson equation domain.}
\label{fig:domainpoisson}
\end{figure}
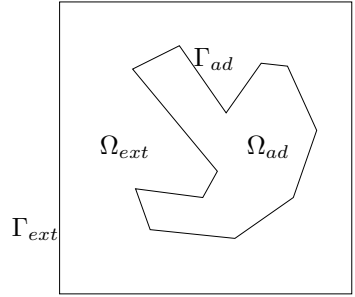

We assume the available sea area to be a simply connected set $\Omega_{ad}$, which we do not require to be convex. 
On the one hand, this allows to consider a larger class of domains of practical interest; on the other hand, it makes the numerical treatment more delicate.\footnote{For example, a convex polygonal set can be described as a system of linear inequalities; this is instead not possible for a concave polygon.} To impose the domain constraint as an equality, we would need a function with value zero inside $\Omega_{ad}$ and positive outside. 
However, such a function would also have zero gradient inside $\Omega_{ad}$, thus making the constraint Jacobian matrix singular and requiring the adoption of a different framework. We instead opt for defining the constraint as an inequality, and requiring the corresponding function $h$ to be negative inside $\Omega_{ad}$. 
This is achieved by solving two Poisson problems (see \cref{fig:domainpoisson}): one on the admissible domain,
\begin{equation}
-\Delta h = \sigma_{int} \, \text{ in $\Omega_{ad}$,} \quad 
h = 0 \, \text{ on $\Gamma_{ad}$},
\label{eq:poissonint}
\end{equation}
and one on an exterior domain $\Omega_{ext}$, 
\begin{equation}
-\Delta h = \sigma_{out}  \,\text{ in $\Omega_{ext}$,} \quad 
h = 0  \,\text{ on $\Gamma_{ad}$,} \quad
\frac{\partial h}{\partial n} = \nu_N  \, \text{ on $\Gamma_{ext}$}.
\label{eq:poissonext}
\end{equation}

We set $\sigma_{int}=-1$, $\sigma_{out}=1$ and $\nu_N=0$.
We now remark on some desirable properties of the solutions.
The maximum principle for both the classical and the weak form of the problems guarantees that $h$ can have no global maxima in the interior of $\Omega_{ad}$ and no global minima in the interior of $\Omega_{ext}$ \cite{Brezis2011}. 
Moreover, for classical solutions it is possible to prove, using a modification of the mean value property \cite{Salsa2016}, that no local maxima can exist in $\Omega_{ad}$ (resp. minima in $\Omega_{ext}$). 
It is also possible to prove, both for the classical and weak formulation, that $\nabla h$ can be zero on at most a set of zero measure.
To compute $h$, the problem needs to be discretized.
If the finite element method is used, with $\mathbb{P}^1$ elements on a Delaunay triangulation, then the numerical solution satisfies local and global discrete maximum properties \cite{Barrenechea_2024}.

The gradient of $h$ needs to be computed to build the differential of constraint functions. If $\mathbb{P}^1$ elements are used, then the exact gradient is discontinuous, potentially making the gradient-flow ODE system stiff. To avoid this possible issue, we can compute a more regular approximation $G_h$ of $\nabla h$ through projection, i.e. by solving the following problem:
\begin{equation*}
(G_h, v) + \eta(\nabla G_h, \nabla v) = (\nabla h, v) \quad \forall v \in H^1(\Omega), \, \Omega=\Omega_{ad} \cup \Omega_{ext}.
\end{equation*}
The constraint of admissible domain is then imposed as
\begin{equation*}
h_{\text{ad}, \ell} = h(\bm{x}_\ell) \neq 0, \quad \ell=1,\dots,N_b.
\end{equation*}

The constraint of minimum distance between devices is written as
\begin{equation}
h_{md, \ell m} = d_{min}^2 - \|\bm{x}_\ell - \bm{x}_m\|^2 \leq 0, \quad \ell, m=1, \dots, N_b, \; l \neq m.
\label{eq:mdconstr}
\end{equation}
It can be verified that it corresponds to a non-convex admissible set \cite{Gambarini2024}.

\subsubsection{Slack variable formulation}
\label{sec:slackwec}
In order to write the final form of our optimal design problem,
we first collect all inequality constraints into a single vector $\bm{h} = [ \bm{h}_{sl}, \bm{h}_{ad}, \bm{h}_{md}]$ and define a corresponding vector of slack variables $\bm{s} = [\bm{s}_{sl}, \bm{s}_{ad}, \bm{s}_{md}]$. The constraints of slamming, admissible domain and minimum distance are then rewritten as
\begin{equation*}
\bm{g}_\mathcal{I}(\bm{u}, \hat{\bm{\gamma}}, \hat{\bm{\zeta}}, \bm{s}) := \bm{h}(\bm{u}, \hat{\bm{\gamma}}, \hat{\bm{\zeta}}) + \bm{s} \odot \bm{s} = 0.
\end{equation*}
By defining a single constraint vector 
\begin{equation}
\bm{g}(\bm{u}, \hat{\bm{\gamma}}, \hat{\bm{\zeta}}, \bm{s}) = [\bm{e}_1(\bm{u}, \hat{\bm{\gamma}}, \hat{\bm{\zeta}}), \dots, \bm{e}_{N_q}(\bm{u}, \hat{\bm{\gamma}}, \hat{\bm{\zeta}}), \bm{g}_\mathcal{I}(\bm{u}, \hat{\bm{\gamma}}, \hat{\bm{\zeta}}, \bm{s})],
\label{eq:defconstrvec}
\end{equation} 
our optimization problem can be restated as an equality-constrained minimization of the opposite of the power \cref{eq:costfunction}:
\begin{equation}
\min_{\bm{u}} \quad f(\bm{u}, \hat{\bm{\gamma}}, \hat{\bm{\zeta}}) = -\frac{1}{2} \sum_{\ell=1}^{N_b} c_\ell \sum_{q=1}^{N_f} \left( \omega_q |\hat{\zeta}_{\ell q}|\right)^2\quad \text{s.t.} \quad
\bm{g}(\bm{u}, \hat{\bm{\gamma}}, \hat{\bm{\zeta}}, \bm{s}) = 0.
\label{eq:optprobcompact}
\end{equation}
For convenience, we finally define $\bm{w} = [\bm{u}, \hat{\bm{\gamma}}, \hat{\bm{\zeta}}, \bm{s}] \in \mathbb{C}^{N_w}$, so that the generic formulation \cref{eq:eqconstropt} is recovered.

\subsection{Algebraic setting and differentials}\label{sec:derivatives}

In this section, we derive the gradient of the cost function and the differential of the constraint vector for our specific application. 
We call $D_w f(\bm{v})$ the Fréchet differential of $f$ at $\bm{w}$ with increment $\bm{v} = [\bm{v}_u, \bm{v}_\gamma, \bm{v}_\zeta, \bm{v}_s]$. 
A description based on the differential operator, instead of a matrix representation, is convenient for computational purposes. 
Indeed, assembling the matrices of the state problem (see \cref{eq:hydoper}, \cref{eq:dynoper}) would require a large number of inefficient matrix-matrix multiplications, and moreover matrices $B_\ell$, $T_{m\ell}$ are sparse, while their product is dense. In our implementation we instead only perform matrix-vector multiplications, by carefully choosing the order of operations: so, the matrices of the state problem are actually treated as abstract linear operators.

The cost function \cref{eq:costfunction} depends explicitly only on the damping coefficients $\bm{c}$ and on the oscillation amplitudes $\hat{\bm{\zeta}}$. We have
\begin{equation*}
\frac{\partial f}{\partial {c_\ell}} = -\frac{1}{2} \sum_{q=1}^{N_f} \left(\omega_q |\hat{\zeta}_{\ell q}|\right)^2, \quad
D_\zeta f(\bm{v}_\zeta) = - \sum_{\ell_1}^{N_b} c_\ell \sum_{q=1}^{N_f} \omega_q^2 \Re\left[ v_{{\zeta}_{\ell q}}^* \hat{\zeta}_{\ell q} \right],
\end{equation*}
where $^*$ denotes the complex conjugate.
The computed differential is $\mathbb{R}$-linear, but not $\mathbb{C}$-linear (in the nomenclature of \cite[sect. III.2]{Lang1987}), because of the presence of the real part operator. 
A suitable setting for our problem is thus the space of complex vectors over the field of real scalars $(\mathbb{C}^{N_w}, \mathbb{R})$, with inner product $(\bm{u}, \bm{v}) = \Re[\bm{u}^H \bm{v}]$, where $^H$ denotes the hermitian (transpose conjugate), and the gradient with respect to oscillation amplitudes is the corresponding Riesz representative of the differential:
\begin{equation*}
\frac{\partial f}{\partial {\hat{\zeta}_{\ell q}}} = -\omega_q^2 c_\ell \hat{\zeta}_{\ell q}.
\end{equation*}
The differential of the constraint vector $\bm{g}$ is split as \begin{equation*}
D_w\bm{g} = [D_w\bm{e}_1, \dots, D_w\bm{e}_{N_q}, D_w\bm{h}_{sl}, D_w\bm{h}_{ad}, D_w\bm{h}_{md}].
\end{equation*} 
Let us start from the residual of the state problem \cref{eq:stateresidual}, focusing on a single frequency and thus disregarding subscript $q$:
\begin{equation*}
D_w \bm{e}(\bm{v}) =
\begin{bmatrix}
M_{\gamma \gamma}(\bm{u}) \bm{v}_\gamma + M_{\gamma \zeta}(\bm{u}) \bm{v}_\zeta + 
D_u \bm{e}_\gamma(\bm{v}_u)\\ 
M_{\zeta \gamma}(\bm{u}) \bm{v}_\gamma + M_{\zeta \zeta}(\bm{u}) \bm{v}_\zeta + 
D_u \bm{e}_\zeta(\bm{v}_u)
\end{bmatrix}.
\end{equation*}
As already mentioned, the blocks of system \cref{eq:stateeq} depend on the decision variables through coordinate transformation matrices, ambient wave vectors and mechanical impedance: in particular, from the $\ell$-th block-row of \cref{eq:hydoper} we obtain
\begin{equation*}
D_u \bm{e}_{\gamma, \ell}(\bm{v}_u) = \sum_{\substack{m=1\\m\neq \ell}}^{N_b} B_\ell \left( D_u T_{m\ell}(\bm{v}_u)\right)^T (\hat{\bm{\gamma}}_m + \bm{R}_m \hat{\zeta}_m) + \frac{H}{2} B_\ell \left(\frac{\partial \bm{a}_\ell}{\partial x_\ell} v_{x_\ell} + \frac{\partial \bm{a}_\ell}{\partial y_\ell} v_{y_\ell} \right),
\end{equation*}
and from the $\ell$-th row of \cref{eq:dynoper}
\begin{equation*}
\begin{split}
D_u e_{\zeta, \ell}(\bm{v}_u) &= \sum_{\substack{m=1\\m\neq \ell}}^{N_b}   \frac{1}{W_\ell} \left[(\widetilde{\bm{Y}}_\ell^D)^T \widetilde{B}_\ell \left(D_u T_{m\ell}(\bm{v}_u)\right)^T    \hat{\bm{\gamma}}_m + \bm{R}_m^T D_u T_{m\ell}(\bm{v}_u) \widetilde{B}_\ell^T \widetilde{\bm{Y}}_\ell^D \hat{\zeta}_m \right] \\
&+ \frac{H}{2W_\ell}\left(\frac{\partial \bm{a}_\ell^T}{\partial x_\ell} v_{x_\ell} + \frac{\partial \bm{a}_\ell^T}{\partial y_\ell} v_{y_\ell} \right) \tilde{B}_\ell^T \tilde{\bm{Y}}_\ell^D - \frac{1}{W_\ell} (e_{\zeta,\ell} - \hat{\zeta}_\ell) D_u W_\ell(\bm{v}_u).
\end{split}
\end{equation*}
The differential of the basis transformation matrix reads
\begin{equation*}
D_u T_{m\ell}(\bm{v}_u) = \frac{\partial T_{m\ell}}{\partial x_\ell} v_{x_\ell} + \frac{\partial T_{m\ell}}{\partial y_\ell} v_{y_\ell} + \frac{\partial T_{m\ell}}{\partial x_m} v_{x_m} + \frac{\partial T_{m\ell}}{\partial y_m} v_{y_m},
\end{equation*}
while the differential of the mechanical impedance \cref{eq:mechimped} is
\begin{equation*}
D_u W_\ell (\bm{v}_u) = \frac{\omega}{\rho g} v_{c_\ell} + \frac{i}{\rho g} v_{\kappa_\ell}.
\end{equation*}
We now report the computation of the required derivatives \cite{jacopo}. The expression of basis transformation matrices is \cite{Child2011}
\begin{equation}
(T_{ij})^{nl}_{mm} = \begin{cases}
\frac{J_l(k r)}{H_n(k r)} H_{n-l}(k L_{ij}) e^{i \alpha_{ij} (n-l)}, & m=0, \\[0.5em]
\frac{I_l(k_m r)}{K_n(k_m r)} K_{n-l}(k_m L_{ij}) e^{i\alpha_{ij}(n-l)} (-1)^l, & m \geq 1,
\end{cases}
\label{eq:basistransmat}
\end{equation}
where $L_{ij} = \|\bm{x}_j - \bm{x}_i\|$, $\alpha_{ij}$ is the angle formed by vector $\bm{x}_j - \bm{x}_i$ and the $x$ axis, $n$, $l$ are indices of the progressive modes, all sharing the same wavenumber $k$, and $m$ is the index of the evanescent mode of wavenumber $k_m$. $J_n$ is a Bessel function of the first kind, $I_n$ is a modified Bessel function of the first kind, $H_n$ is a Hankel function of the first kind and $K_n$ is a modified Bessel function of the second kind, subscripts indicating the order.  Differentiating \cref{eq:basistransmat} and using the expression for the derivative of $H_n$ \cite{Arfken2005}:
\begin{equation*}
\frac{dH_n(x)}{dx} = \frac{H_{n-1}(x) - H_{n+1}(x)}{2},
\end{equation*} 
the same holding for $K_n$, yields \small
\begin{equation*}
\left( \frac{\partial T_{ij}}{\partial x_\ell} \right)^{nl}_{mm} = \begin{cases}
\begin{aligned}
\frac{J_l(k r)}{H_n(k r)} &\left[ k \frac{\partial L_{ij}}{\partial x_\ell} \frac{H_{n-l-1}(kL_{ij}) - H_{n-l+1}(kL_{ij})}{2} \right. \\ & \left. + i \frac{\partial \alpha_{ij}}{\partial x_\ell} (n-l) H_{n-l}(k L_{ij})\right] e^{i \alpha_{ij} (n-l)},\end{aligned} & m=0, \\
\begin{aligned}
\frac{I_l(k_m r)}{K_n(k_m r)} &\left[k_m \frac{\partial L_{ij}}{\partial x_\ell} \frac{K_{n-l-1}(kL_{ij}) - K_{n-l+1}(kL_{ij})}{2} \right. \\ &
\left. + i \frac{\alpha_{ij}}{\partial x_\ell} K_{n-l}(k_m L_{ij}) \right] e^{i \alpha_{ij}(n-l)}(-1)^l, 
\end{aligned} & m \geq 1 .
\end{cases}
\end{equation*} \normalsize
For the derivatives of $L_{ij}$ and $\alpha_{ij}$, we start from $
L_{ij} \cos \alpha_{ij} = x_j - x_i$, $L_{ij} \sin \alpha_{ij} = y_j - y_i$.
Differentiation with respect to $x_i$ yields the system
\begin{equation*}
\begin{bmatrix}
\cos \alpha_{ij} & -L_{ij} \sin \alpha_{ij} \\
\sin \alpha_{ij} & L_{ij} \cos \alpha_{ij}
\end{bmatrix}
\begin{bmatrix}
\dfrac{\partial L_{ij}}{\partial x_i} \\[0.8em]
\dfrac{\partial \alpha_{ij}}{\partial x_i}
\end{bmatrix}
=
\begin{bmatrix}
-1 \\ 0
\end{bmatrix},
\end{equation*}
whose matrix has determinant $L_{ij}$ and thus is always invertible for distinct points.
Analogous systems hold for the other coordinates involved.
The expressions of the ambient wave coefficients are \cite{Child2011}
\begin{equation*}
(\bm{a}_\ell)^n_m = \begin{cases}
J_n(k r) e^{i n (\pi/2 - \beta)} \Theta_\ell & m=0, \\
0, & m\geq 0,
\end{cases} \quad \Theta_\ell=e^{i k (x_\ell \cos\beta + y_\ell \sin \beta)},
\end{equation*}
and their derivatives are
\begin{equation*}
\dfrac{\partial \bm{a}_\ell}{\partial x_\ell} = i k \cos\beta \bm{a}_\ell, \quad
\dfrac{\partial \bm{a}_\ell}{\partial y_\ell} = i k \sin\beta \bm{a}_\ell.
\end{equation*} 
The slamming constraint vector $\bm{h}_{sl}$ depends on the oscillation amplitudes $\hat{\bm{\zeta}}$ explicitly, and on the devices' positions through the phase of incident waves. Its derivatives with respect to positions can be directly evaluated (see \cref{eq:etadecomp}):
\begin{align*}
\frac{\partial h_{sl,\ell}}{\partial x_\ell} &= -2 \sum_{q=1}^{N_f} \Re\left[ (\hat{\zeta}_{q\ell} - \hat{\eta}_{q\ell})^* \frac{\partial \hat{\eta}_{q\ell}}{\partial x_\ell} \right], \\
\frac{\partial \hat{\eta}_{q\ell}}{	\partial x_\ell} &= \frac{\partial}{\partial x_\ell} \left[ i \frac{H_q}{2} \exp[i k (x_\ell \cos\beta + y_\ell \sin \beta)]\right] = i k \cos\beta \hat{\eta}_{q\ell}.
\end{align*}
We observe that 
\begin{equation*}
\Re\left[ \hat{\eta}_{q\ell}^* \frac{\partial \hat{\eta}_{q\ell}}{\partial x_\ell} \right] = \frac{1}{2} \frac{\partial}{\partial x_\ell} \left|\hat{\eta}_{q\ell} \right|^2 = 0,
\end{equation*}
since coordinate $x_\ell$ has an effect on the phase of $\hat{\eta}_{q\ell}$, but not on its magnitude. Hence,
\begin{equation*}
\frac{\partial \hat{\eta}_{q\ell}}{\partial x_\ell} = -2  \cos\beta \sum_{q=1}^{N_f}  k_q \Re \left( i \hat{\zeta}_{q\ell}^* \hat{\eta}_{q\ell} \right) = 2  \cos\beta \sum_{q=1}^{N_f}  k_q \Re \left( i \hat{\zeta}_{q\ell} \hat{\eta}_{q\ell}^* \right).
\end{equation*}
An analogous result holds for $y_\ell$, with $\sin\beta$ in place of $\cos\beta$.
Differentiation with respect to oscillation amplitudes yields
\begin{equation*}
D_\zeta h_{sl,\ell}(\bm{v}_{\zeta}) = \sum_{q=1}^{N_f} 2 \Re \left[ v_{{\zeta}_{q\ell}}^* \left( \hat{\zeta}_{q\ell} - \hat{\eta}_{q\ell} \right) \right].
\end{equation*}
The full differential of the $\ell$-th component of the slamming constraint reads
\begin{equation*}
D_w h_{sl, \ell}(\bm{v}) = \frac{\partial h_{sl, \ell}}{\partial x_\ell} v_{x_\ell} + \frac{\partial h_{sl, \ell}}{\partial y_\ell} v_{y_\ell} + \sum_{q=1}^{N_f} 2 \Re \left[ v_{{\zeta}_{q\ell}}^* \left( \hat{\zeta}_{q\ell} - \hat{\eta}_{q\ell} \right) \right] + 2 s_{sl, \ell} v_{s_{sl,\ell}}.
\end{equation*}
The computation of the gradient of the admissible domain constraint has been discussed in \cref{sec:domconstr}, and the corresponding differential is
\begin{equation*}
D_w h_{ad, \ell}(\bm{v}) = G_h(\bm{x}_\ell) \cdot \bm{v}_{x_\ell} + 2 s_{ad, \ell} v_{s_{ad, \ell}}. 
\end{equation*}
The derivatives of the minimum distance constraint function \cref{eq:mdconstr} with respect to the coordinates can be readily computed, yielding the differential
\begin{equation*}
D_w h_{md, \ell m}(\bm{v}) = -2 (\bm{x}_\ell - \bm{x}_m) \cdot (\bm{v}_{x_\ell} - \bm{v}_{x_m}) + 2 s_{md, \ell m} v_{s_{md, \ell m}}.
\end{equation*}
The computation of the adjoint of the differential, $(D_w\bm{g})'$, appearing in \cref{eq:gradflowtanabe}, can be done directly using the definition of adjoint operator: $(\bm{p}, D_w\bm{g}(\bm{v})) = (\bm{v}, (D_w\bm{g})'(\bm{p}))$. We define $\bm{p} = [\bm{p}_e, \bm{p}_{sl}, \bm{p}_{ad}, \bm{p}_{md}]$, and hence the computation amounts to writing $\Re[\bm{p}^H D_w\bm{g}(\bm{v})]$, taking the hermitian inside the real part operator and collecting terms. 
As for $D_w\bm{g}$, we split the computation of the adjoint operator into block-rows, corresponding to the decision, state and slack variables.
\begin{align*}
(D_w \bm{g})'_{x_\ell} (\bm{p}) &= \sum_{q=1}^{N_f} \Re \biggl\{ \biggl[ \sum_{\substack{m=1\\m\neq \ell}}^{N_b} \left( \hat{\bm{\gamma}}_m^H + \hat{\zeta}_m^* \bm{R}_m^H \right) \frac{\partial T_{m\ell}^*}{\partial x_\ell} B_\ell^H + \frac{H}{2} \frac{\partial \bm{a}_\ell^H}{\partial x_\ell} B_\ell^H \biggr] \bm{p}_{e_{\gamma,\ell}} \biggr. \nonumber\\
&+ \biggl[ \sum_{\substack{m=1\\m\neq \ell}}^{N_b} \left( \hat{\bm{\gamma}}_\ell^H + \hat{\zeta}_\ell^* \bm{R}_\ell^H \right) \frac{\partial T_{\ell m}^*}{\partial x_\ell} B_m^H  \biggr] \bm{p}_{e_{\gamma, m}} + \frac{1}{W_\ell}\biggl[ \sum_{\substack{m=1\\m\neq \ell}}^{N_b}  \hat{\bm{\gamma}}_m^H \frac{\partial T_{m\ell}^*}{\partial x_\ell} \tilde{B}_\ell^H (\tilde{\bm{Y}}_\ell^D)^* \biggr. \nonumber\\&\biggl. +  \sum_{\substack{m=1\\m\neq \ell}}^{N_b}  \hat{\zeta}_m^* (\bm{Y}_\ell^D)^H \tilde{B}_\ell^* \frac{\partial T_{m\ell}^H}{\partial x_\ell} \bm{R}_m^* \biggr] p_{e_{\zeta,\ell}} + \frac{H}{2W_\ell} (\tilde{\bm{Y}}_\ell^D)^H \tilde{B}_\ell^* \frac{\partial \bm{a}_\ell^*}{\partial x_\ell} p_{e_{\zeta,\ell}} \\
&\biggl. +\frac{1}{W_m} \biggl[ \sum_{\substack{m=1\\m\neq \ell}}^{N_b}  \hat{\bm{\gamma}}_\ell^H \frac{\partial T_{\ell m}^*}{\partial x_\ell} \tilde{B}_m^H (\tilde{\bm{Y}}_m^D)^* +  \sum_{\substack{m=1\\m\neq \ell}}^{N_b}  \hat{\zeta}_\ell^* (\bm{Y}_m^D)^H \tilde{B}_m^* \frac{\partial T_{\ell m}^H}{\partial x_\ell} \bm{R}_\ell^* \biggr] p_{e_{\zeta,m}} \biggr\}_q \nonumber\\
& + \frac{\partial h_{sl,\ell}}{\partial x_\ell} p_{sl, \ell} + G_h(\bm{x}_\ell)\cdot \hat{\bm{x}} p_{ad, \ell} - 2 \sum_{\substack{m=1\\m\neq \ell}}^{N_b} (x_\ell - x_m) p_{md, \ell m} \nonumber ,
\end{align*}
where $\hat{\bm{x}}$ denotes the unit vector in the $x$ direction, and where all the quantities inside curly brackets are related to the $q$-th frequency component, as indicated by the subscript at the end. The same relation holds for $(D_u \bm{g})'_{y_\ell}$, replacing $x$ with $y$. The component related to damping is
\begin{equation*}
(D_w \bm{g})'_{c_\ell} (\bm{p}) = \sum_{q=1}^{N_f} \Re \left[ -\frac{1}{W_\ell} (e_{\zeta, \ell} - \hat{\zeta}_\ell)^* \frac{\partial W_\ell^*}{\partial c_\ell} p_{e\zeta, \ell} \right]_q,
\end{equation*}
and the one related to stiffness is analogous, with $\kappa_\ell$ in place of $c_\ell$. The rows corresponding to state variables $\hat{\bm{\gamma}}$, $\hat{\bm{\zeta}}$ are
\begin{align*}
(D_w\bm{g})'_{\gamma}(\bm{p}) &= \sum_{q=1}^{N_f} \left[ M_{\gamma\gamma}^H \bm{p}_{e\gamma} + M_{\zeta\gamma}^H \bm{p}_{e\zeta}\right]_q \\
(D_w\bm{g})'_{\zeta}(\bm{p}) &= \sum_{q=1}^{N_f} \left[ M_{\gamma\zeta}^H \bm{p}_{e\gamma} + M_{\zeta\zeta}^H \bm{p}_{e\zeta} \right]_q + 2 \sum_{\ell=1}^{N_b} \sum_{q=1}^{N_f} (\hat{\zeta}_{q\ell} - \hat{\eta}_{q\ell}) p_{sl, \ell},
\end{align*}
and finally the ones corresponding to slack variables are
\begin{equation*}
(D_w\bm{g})'_{s_{sl}} = 2 \bm{s}_{sl} \odot \bm{p}_{sl},
\quad(D_w\bm{g})'_{s_{ad}} = 2 \bm{s}_{ad} \odot \bm{p}_{ad}, \quad(D_w\bm{g})'_{s_{md}} = 2 \bm{s}_{md} \odot \bm{p}_{md}.
\end{equation*}

Some remarks are in order regarding the scaling described in \cref{sec:scaling}.
Rows of $J_g$ corresponding to inequality constraints are normalized with their $2$-norm, while block-rows corresponding to the state problem are left unchanged for two reasons. The first is that the state problem is already fairly well conditioned \cite[Sec. 4.3]{Child2011}; the second is that, as mentioned in \cref{sec:derivatives}, the blocks of the state system are implemented as linear operators instead of matrices for computational efficiency, so that the rows of their representatives should be recovered by computing their action on all vectors of the canonical basis, an expensive step.


\section{Numerical experiments}\label{sec:numtests}
\begin{figure}[h!]
\centering
\subfloat{
\begin{overpic}[scale=.25]{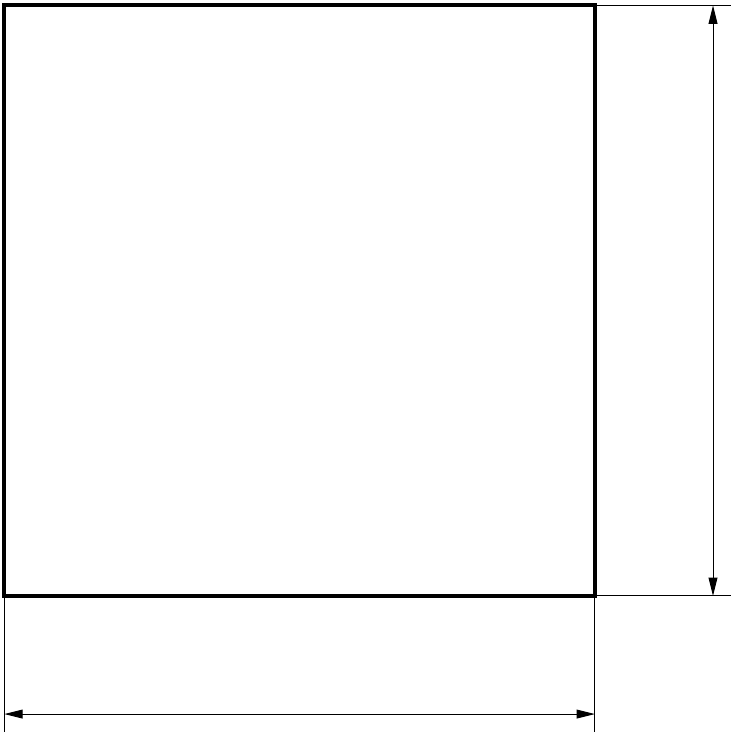}
\footnotesize
\put(90,55){\rotatebox{90}{50}}
\put(37, 5){50}
\end{overpic}
}
\hspace{0.1\linewidth}
\subfloat{
\begin{overpic}[scale=.25]{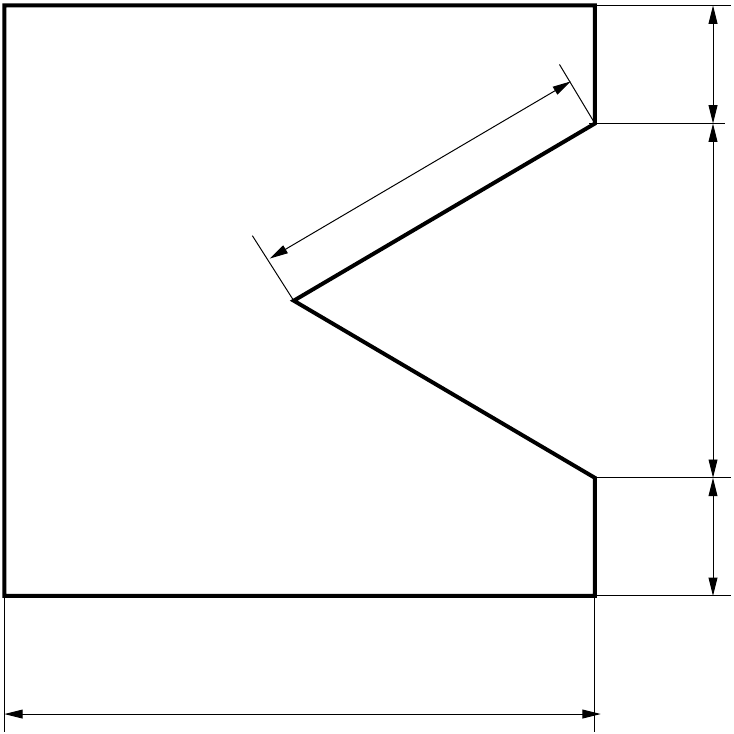}
\footnotesize
\put(37, 5){50}
\put(90,87){\rotatebox{90}{10}}
\put(90,55){\rotatebox{90}{30}}
\put(90,23){\rotatebox{90}{10}}
\put(51,76){\rotatebox{30}{30}}
\end{overpic}
}
\caption{Available sea area for the square (left) and cut square (right) test cases; all lengths in meters.}
\label{fig:availableareas}
\end{figure}
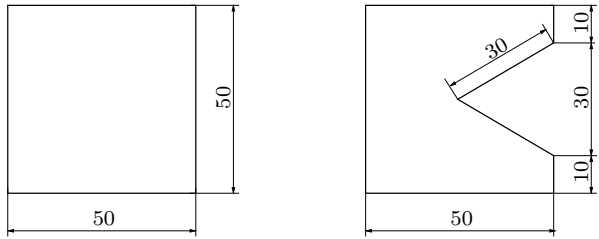

In this section, we present numerical experiments based on two admissible domains, represented in \cref{fig:availableareas}. The first is a simple square, while the second is a square with a triangular cut. The latter represents an available sea area in which a floating wind turbine platform, whose typical shape is a triangle (we mention as an example the Windfloat project \cite{Roddier_2010}), needs to be installed, thus reducing the space available for the wave energy park. The resulting domain is not convex. 

To compare different possible strategies for the numerical solution of the optimization problem, four simulations were performed for each domain geometry, using the following settings:
\begin{itemize}
\item [S1:] Explicit Euler, $\Delta t = 1$, $\tau_{CG} = 10^{-6}$;
\item [S2:] Explicit Euler, $\Delta t = 1.5$, $\tau_{CG} = 10^{-6}$;
\item [S3:] RK12, adaptive $\Delta t$: $\tau_{RK}^{rel}= 10^{-3}$, $\tau_{RK}^{abs}= 10^{-6}$, $\tau_{CG} = 10^{-6}$;
\item [S4:] RK12, adaptive tolerances for RK and CG.
\end{itemize}
In all computations we consider devices of radius $R = 2 \text{ m}$ and draft $d = 0.5 \text{ m}$, and a sea depth $D = 30 \text{ m}$. The sea state is represented by a Pierson-Moskowitz wave spectrum \cite{Falnes2020} with $T_e=8\text{ s}$, $H_s = 2.12\text{ m}$ and waves directed along the positive $x$ coordinate ($\beta = 0$). 
The spectrum is discretized with 30 harmonics. 
For the solution of the state problem, $N_n=4$, $N_m=25$ are chosen.
The stopping tolerance is set to $\tau_\Psi = 10^{-3}$.

Numerical simulations have been performed on a laptop with a 4-core, 8-thread Intel Core i5-10210U CPU and 8 GB of RAM. The presented test cases require about 2 GB of RAM for execution. The code has been implemented in Python and has been made available for reproducibility of the results at \url{https://github.com/marcogambarini/floatcyl}.

\begin{table}[]
\centering
\begin{tabular}{lccccc}
& $f_{end}$ & $\|g\|_{end}$ & $\|\bm{\Psi}\|_{end}$ & ncalls & tot Tcalls (s) \\\cline{2-6}\noalign{\vspace{0.25ex}}
run1 & \sisetup{round-mode = places, round-precision = 3}
\num{-1.3877642275484887} & \sisetup{round-mode = places, round-precision = 3}
\num{1.850645839574247e-06} & \sisetup{round-mode = places, round-precision = 4}
\num{0.0009832888503667825} & 406 & \sisetup{round-mode = places, round-precision = 0}
\num{11251.46498862699}\\
run2 & \sisetup{round-mode = places, round-precision = 3}
\num{-1.3877652265303537} & \sisetup{round-mode = places, round-precision = 3}
\num{1.684830976322027e-06} & \sisetup{round-mode = places, round-precision = 4}
\num{0.0009711906094975896} & 269 & \sisetup{round-mode = places, round-precision = 0}
\num{7672.856474587974}\\
run3 & \sisetup{round-mode = places, round-precision = 3}
\num{-1.3877600527041514} & \sisetup{round-mode = places, round-precision = 3}
\num{2.610744634455896e-06} & \sisetup{round-mode = places, round-precision = 4}
\num{0.0009799064878732707} & 252 & \sisetup{round-mode = places, round-precision = 0}
\num{7272.766847484985}\\
run4 & \sisetup{round-mode = places, round-precision = 3}
\num{-1.387765094697356} & \sisetup{round-mode = places, round-precision = 3}
\num{1.6395578890907083e-06} & \sisetup{round-mode = places, round-precision = 4}
\num{0.0009549837076730104} & 242 & \sisetup{round-mode = places, round-precision = 0}
\num{6161.8742464309835}\\
\end{tabular}
\caption{Results for the square admissible domain}
\label{tab:res-square}
\end{table}

\begin{table}[]
\centering
\begin{tabular}{lccccc}
& $f_{end}$ & $\|g\|_{end}$ & $\|\bm{\Psi}\|_{end}$ & ncalls & tot Tcalls (s) \\\cline{2-6}\noalign{\vspace{0.25ex}}
run1 & \sisetup{round-mode = places, round-precision = 3}
\num{-1.3956946350751818} & \sisetup{round-mode = places, round-precision = 3}
\num{2.716152002074453e-06} & \sisetup{round-mode = places, round-precision = 4}
\num{0.0009951476086842836} & 968 & \sisetup{round-mode = places, round-precision = 0}
\num{25154.509994049018}\\
run2 & \sisetup{round-mode = places, round-precision = 3}
\num{-1.3952050676578278} & \sisetup{round-mode = places, round-precision = 3}
\num{4.033775003397032e-06} & \sisetup{round-mode = places, round-precision = 4}
\num{0.0009875682018874832} & 673 & \sisetup{round-mode = places, round-precision = 0}
\num{18306.971098903916}\\
run3 & \sisetup{round-mode = places, round-precision = 3}
\num{-1.3953589009875065} & \sisetup{round-mode = places, round-precision = 3}
\num{0.0009349979252956905} & \sisetup{round-mode = places, round-precision = 4}
{\color{red}\num{0.003833965314450835}} & 1112 & \sisetup{round-mode = places, round-precision = 0}
\num{31498.80963385893}\\
run4 & \sisetup{round-mode = places, round-precision = 3}
\num{-1.3951070216194237} & \sisetup{round-mode = places, round-precision = 3}
\num{4.4987728669244885e-05} & \sisetup{round-mode = places, round-precision = 4}
\num{0.0009950476034428688} & 763 & \sisetup{round-mode = places, round-precision = 0}
\num{15019.599370843018}\\
\end{tabular}
\caption{Results for the cut square test case}
\label{tab:res-cutsquare}
\end{table}

The results of the test cases are reported in \cref{tab:res-square} and \cref{tab:res-cutsquare}: we show the final values of cost, $f_{end}$, constraint vector norm, $\|\bm{g}\|_{end}$, and stopping indicator, $\|\bm{\Psi}\|_{end}$. Values of $\|\bm{\Psi}\|_{end}$ in red signify that the stopping criterion has not been met within the maximal fictitious time interval. We further report the total number of calls to the function computing $\bm{\Psi}$ and the total time required by such calls.
For both cases, using the explicit Euler method (S1 and S2) leads to convergence. 
The Runge-Kutta method with adaptive tolerance (S4) also reaches convergence, requiring a significantly smaller computational time. The Runge-Kutta method without adaptive tolerances (S3), instead, is only able to fulfill the stopping criterion in the square case. Regarding the cost function, we observe that the result reached for the cut square geometry is more favorable than the one on the square. 
This may appear counterintuitive; for this reason, a test using the result on the cut square geometry as initial guess for an optimization on the whole square was performed. The resulting optimal solution has a negligible difference in terms of cost from the initial guess, suggesting that the solution obtained in the cut square case is also close to a local minimum over the entire square.

\begin{figure}[]
\centering
\subfloat{
\includegraphics[width=0.47\textwidth]{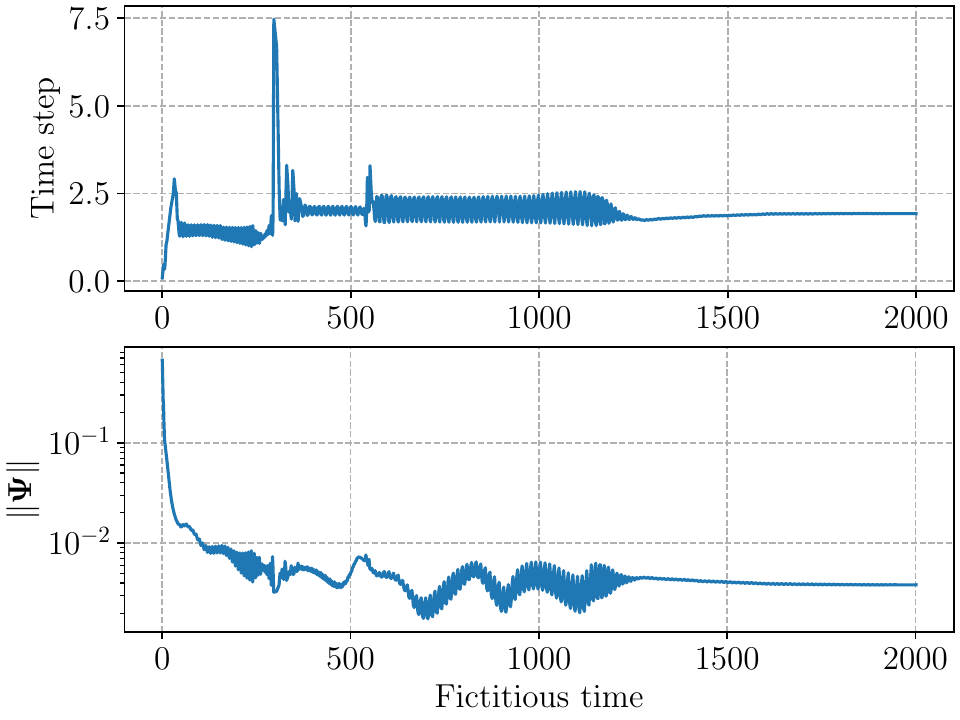}}
\hfill
\subfloat{
\includegraphics[width=0.47\textwidth]{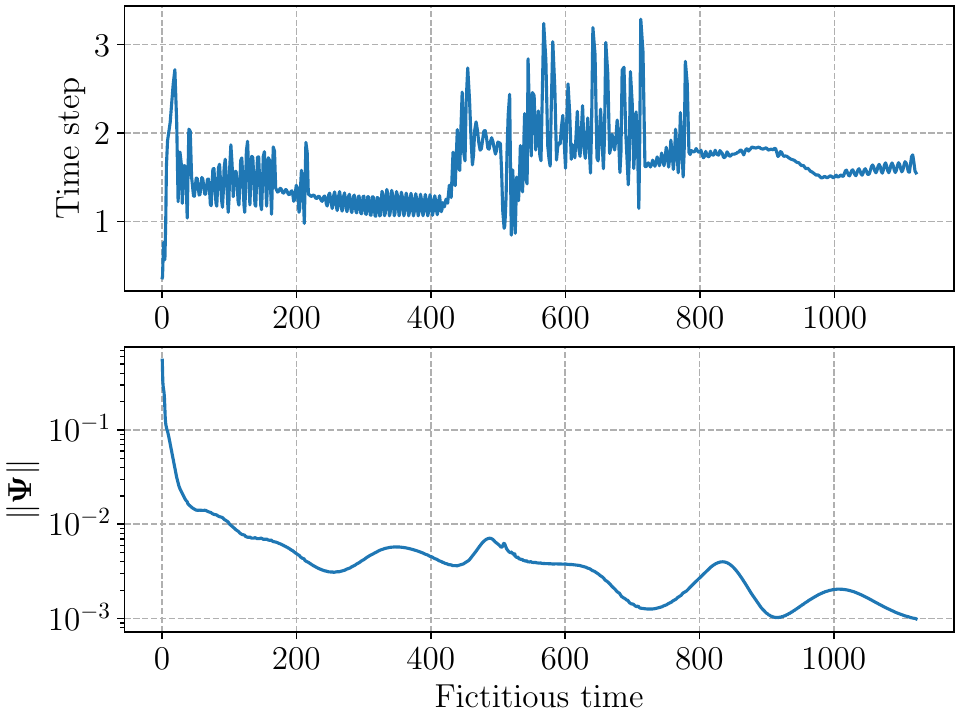}}
\caption{Cut square geometry, history of adaptive time step and stopping indicator, settings S3 (left) and S4 (right)}\label{fig:cutsquare-test34-dtphi}
\end{figure}

\begin{figure}[]
\centering
\subfloat{
\includegraphics[width=0.47\textwidth]{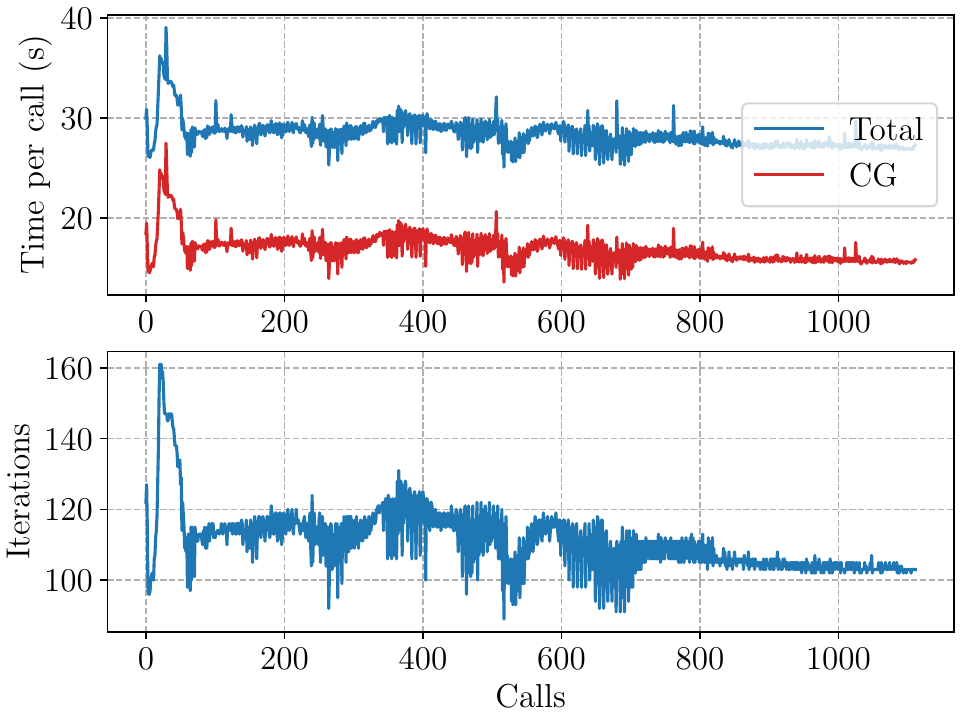}}
\hfill
\subfloat{
\includegraphics[width=0.47\textwidth]{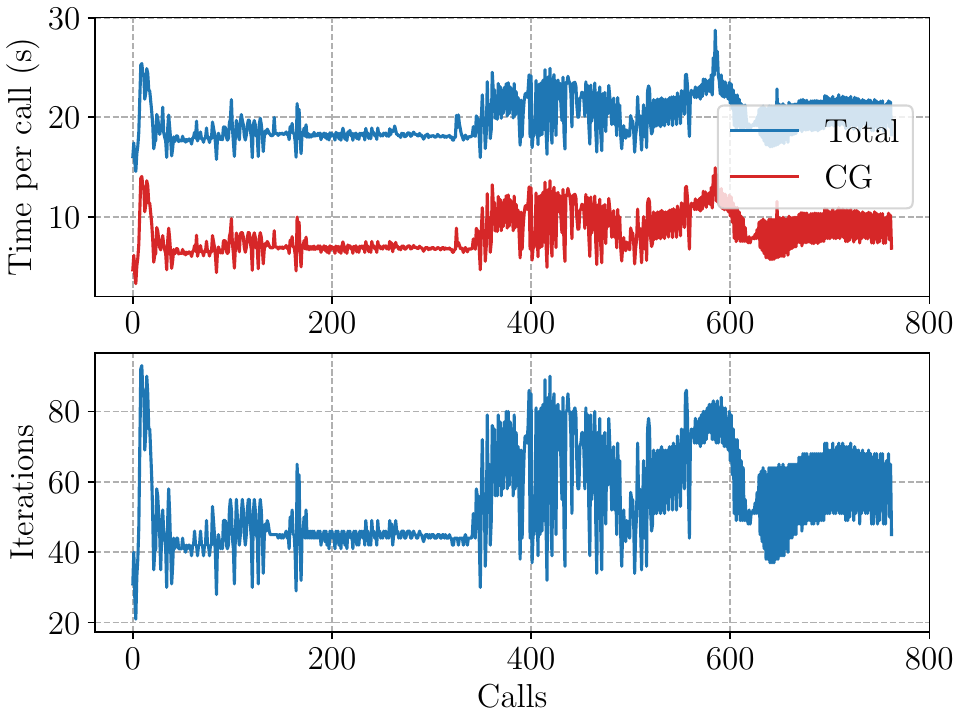}}
\caption{Cut square geometry, history of walltime per call and CG iterations, settings S3 (left) and S4 (right)}\label{fig:cutsquare-test34-monitor}
\end{figure}

We now focus on the Runge-Kutta method applied to the cut square test case. \cref{fig:cutsquare-test34-dtphi} shows the evolution of adaptive time step and stopping indicator. 
For both settings S3 and S4, the selected time step oscillates around 2 for most of the time, with more pronounced reductions for S4 in the first half of the simulation. The value of $\|\bm{\Psi}\|$ is not monotonically decreasing, nor is it expected to be. However, the adaptive tolerance criteria adopted in S4 enable converging to the prescribed stopping tolerance, while in the case of S3 the stopping indicator, after an oscillatory phase, decreases very slowly and seems to stagnate.
We further show, in \cref{fig:cutsquare-test34-monitor}, that using the adaptive tolerance strategy (settings S4) significantly reduces the number of iterations of the conjugate gradient method compared to using fixed tolerances (settings S3). The time required for computing $\bm{\Psi}$ at each time step is also reduced, as expected.

\begin{figure}
\centering
\subfloat{
\includegraphics[width=0.32\textwidth, trim={1cm 0. 2.5cm 0}, clip]{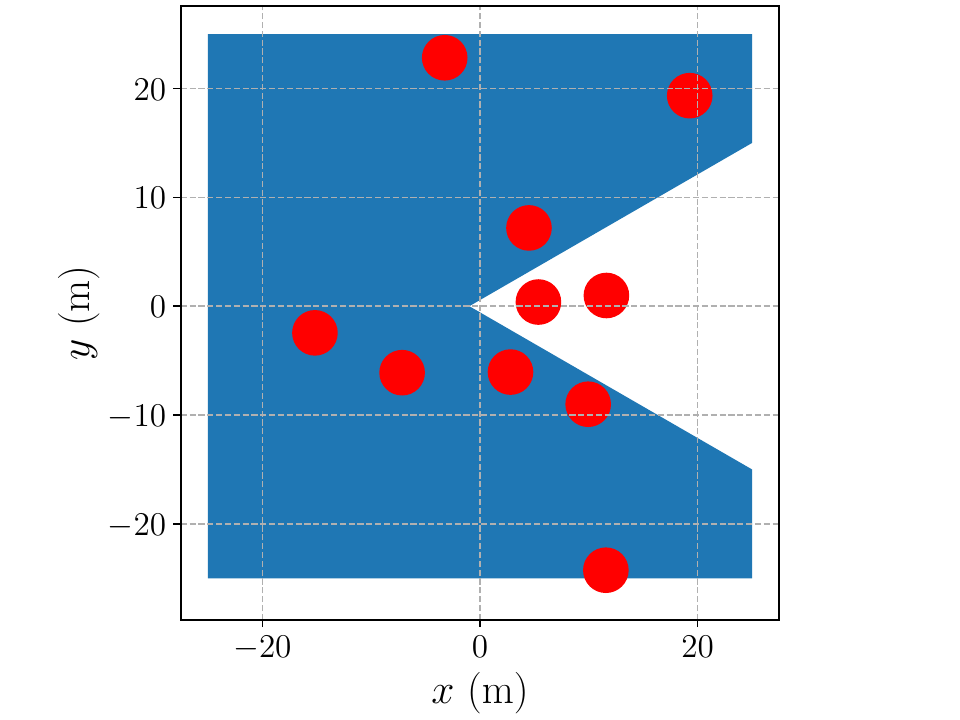}}
\hfill
\subfloat{
\includegraphics[width=0.32\textwidth, trim={1cm 0. 2.5cm 0}, clip]{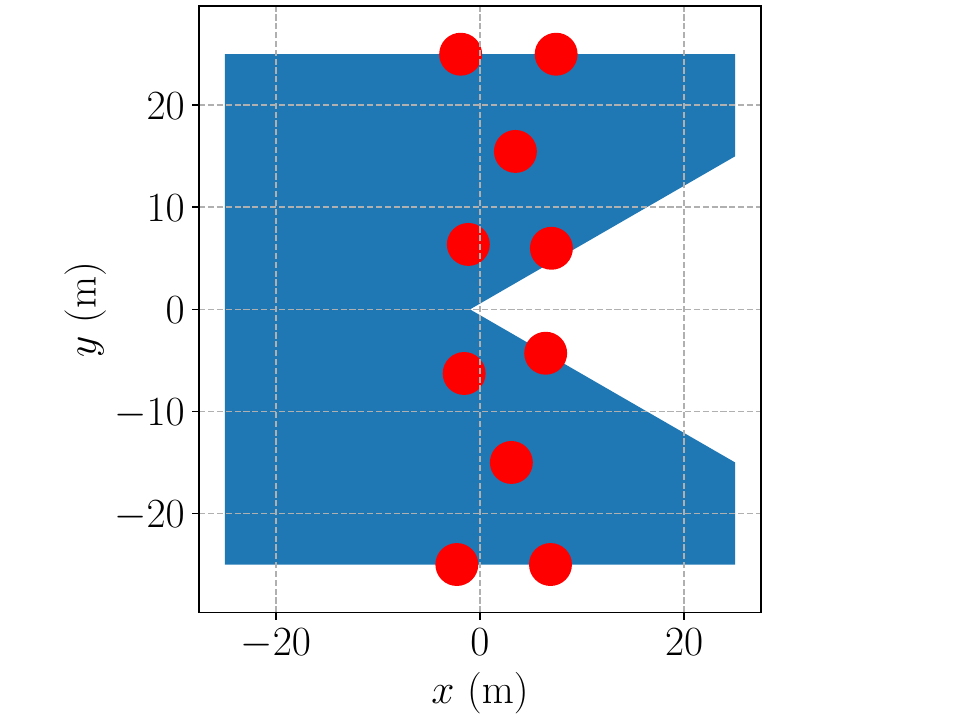}}
\hfill
\subfloat{
\includegraphics[width=0.32\textwidth, trim={1cm 0. 2.5cm 0}, clip]{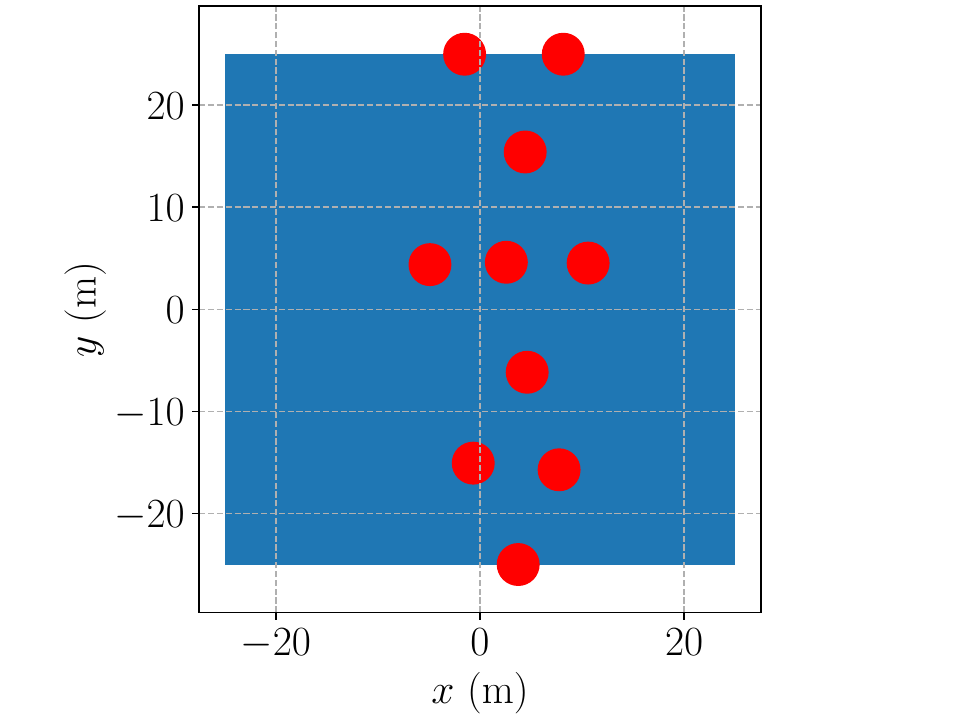}}
\caption{Initial layout guess (left), final layout for the cut square geometry (center), and for the square geometry (right), with settings S4}\label{fig:test4-layout}
\end{figure}

\begin{figure}[]
\centering
\subfloat{
\includegraphics[trim={3cm, 0, 0, 0}, clip, width=0.46\textwidth]{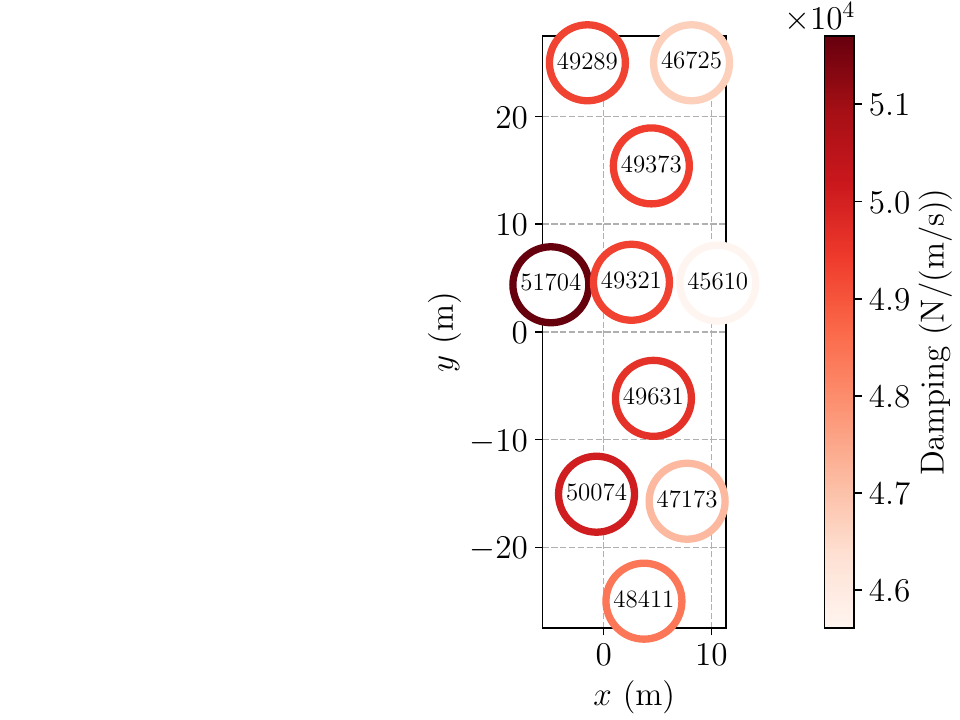}}
\hfill
\subfloat{
\includegraphics[trim={3cm, 0, 0, 0}, clip, width=0.45\textwidth]{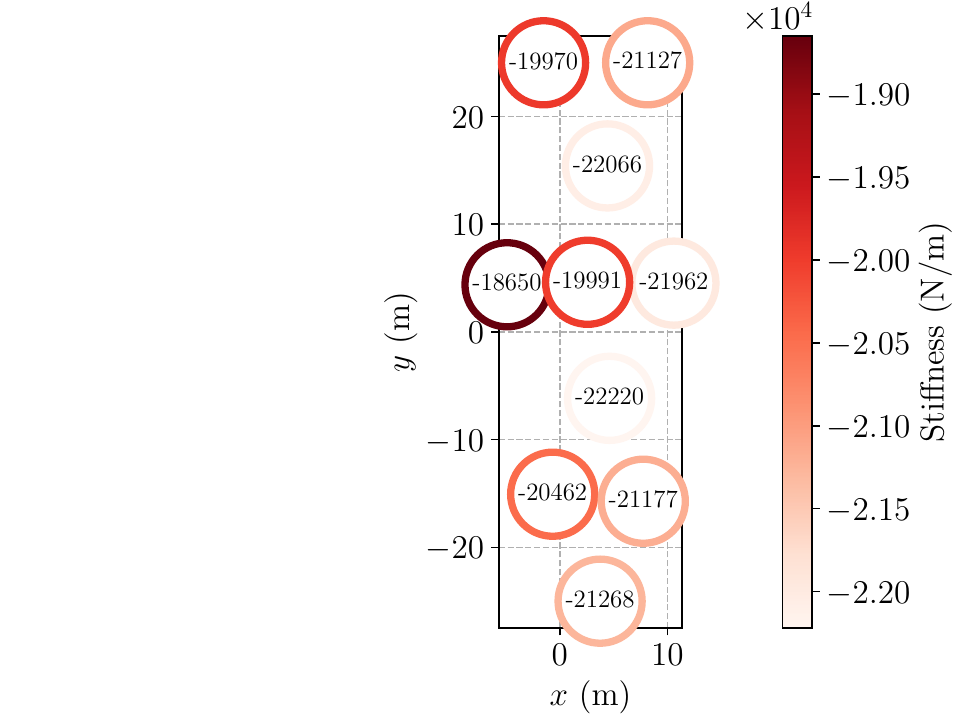}}
\caption{Square geometry, settings S4, optimal damping (left) and stiffness (right)}\label{fig:square-test4-controls}
\end{figure}

\begin{figure}
\centering
\subfloat{
\includegraphics[trim={3cm, 0, 0, 0}, clip, width=0.45\textwidth]{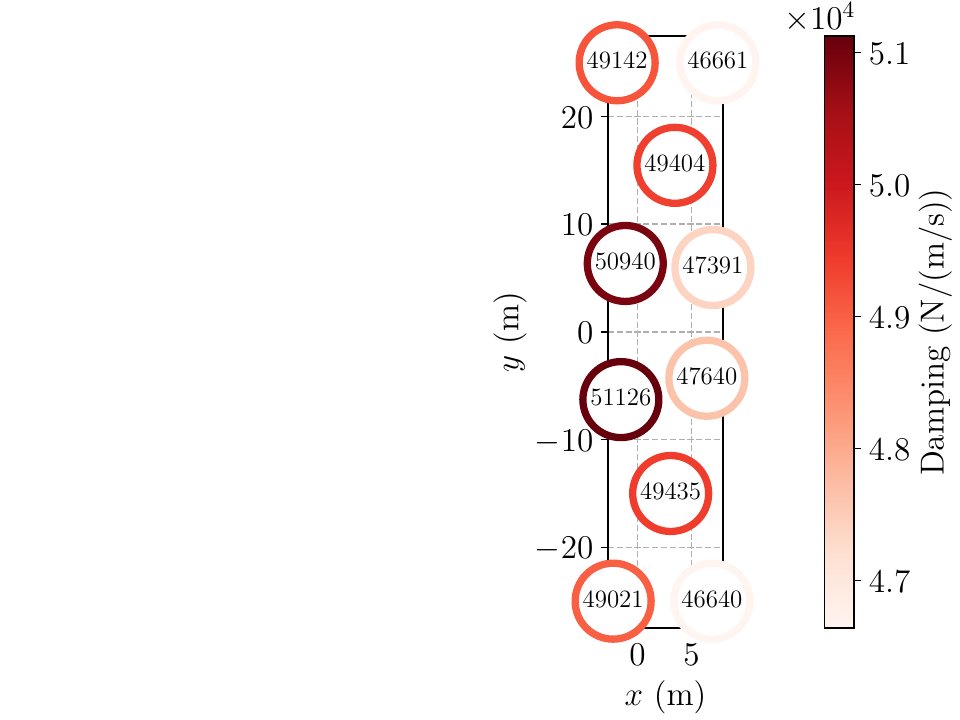}}
\hfill
\subfloat{
\includegraphics[trim={3cm, 0, 0, 0}, clip, width=0.45\textwidth]{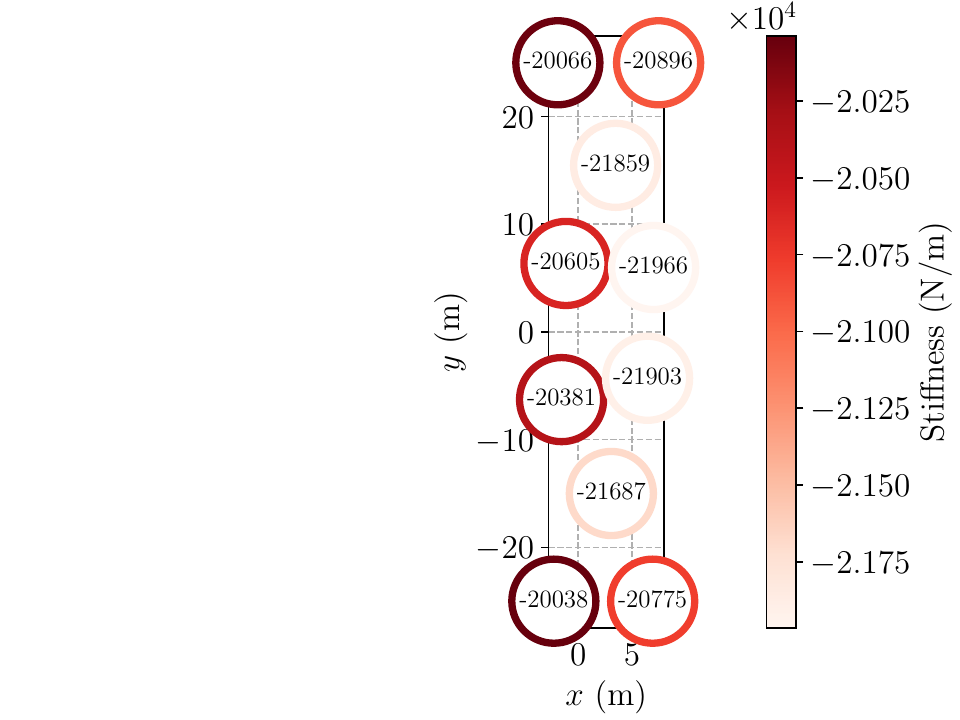}}
\caption{Cut square geometry, settings S4, optimal damping (left) and stiffness (right)}\label{fig:cutsquare-test4-controls}
\end{figure}

The results of optimization in terms of positions are reported in \cref{fig:test4-layout}, while results for controls are shown in  \cref{fig:square-test4-controls} for the square case, and in \cref{fig:cutsquare-test4-controls} for the cut square test case.
The results show that the treatment of the domain constraint is effective even when the initial guess is not feasible. 
Moreover, upwave devices are assigned greater values of both stiffness and damping than downwave devices.
A physical interpretation of this fact is that upwave devices should be tuned to higher frequencies than downwave devices.
The same qualitative result was obtained in \cite{Gambarini2023}.


\begin{figure}[h!]
    \centering
    \subfloat{
    \includegraphics[width=0.48\linewidth]{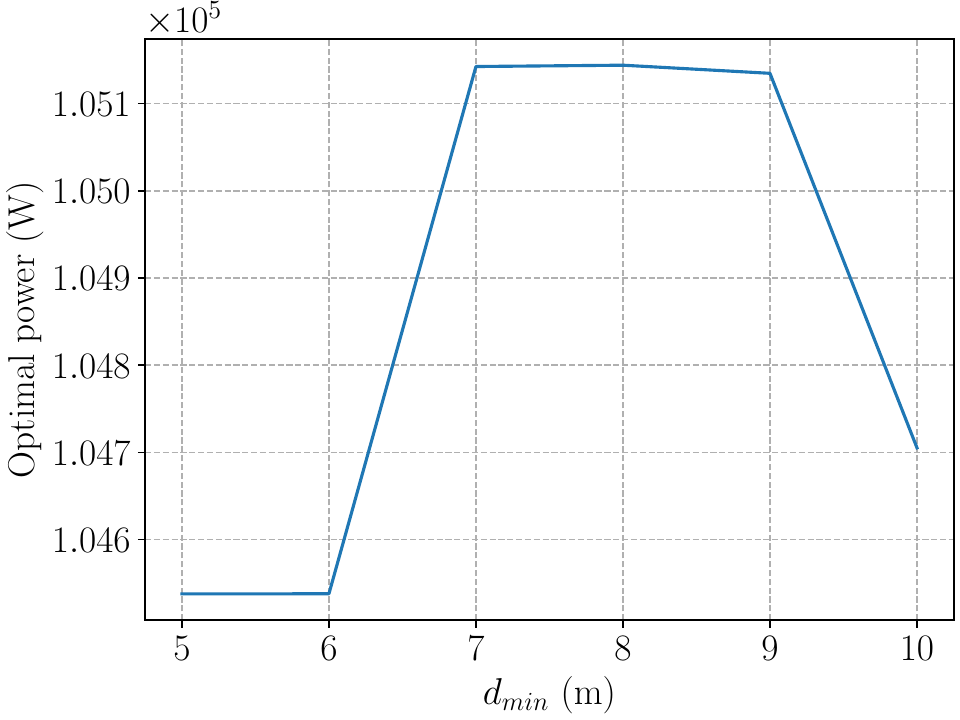}}
    \hfill
    \subfloat{
    \includegraphics[width=0.48\linewidth]{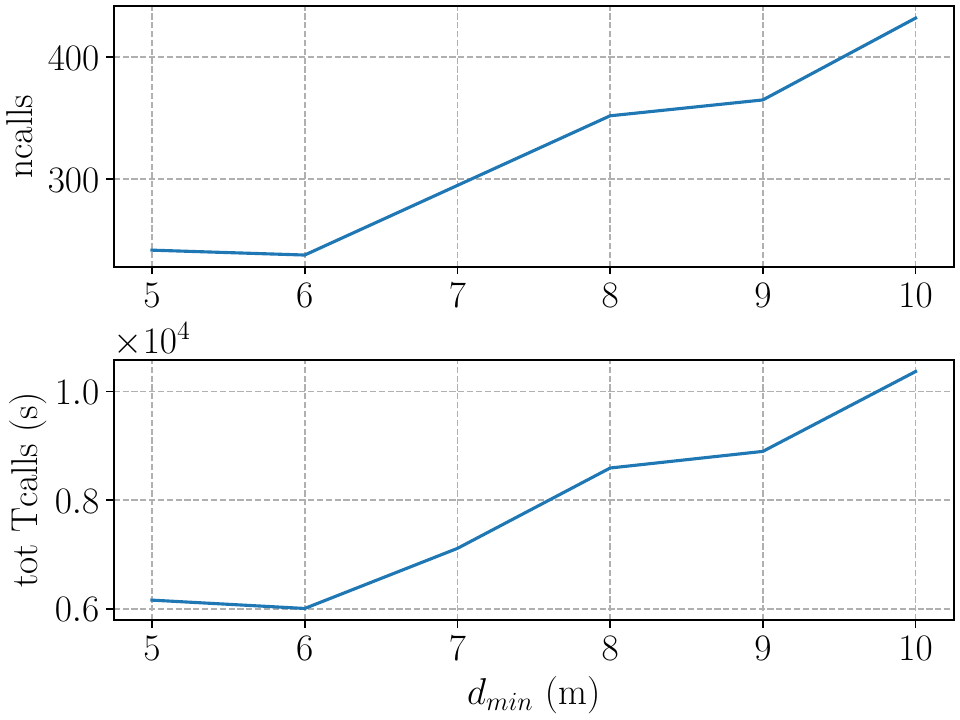}}
    \caption{Square geometry, 10 devices, variable value of $d_{min}$}
    \label{fig:vardist}
\end{figure}

From the point of view of robustness with respect to changes in the constraint parameters, we show in \cref{fig:vardist} the results obtained by varying the minimum distance between devices, $d_{min}$. The initial layouts are compliant with the constraint, so that they are all different. Correspondingly, different results in terms of layout and controls are obtained; however, as the figure shows, the variation in optimal power is very small. In terms of computational cost, the number of calls to the computation of $\bm{\Psi}$ increases as the constraint is tightened, and the computational time grows approximately proportionally, but the proposed method is still able to reach convergence up to the required tolerance.

\section{Conclusions}\label{sec:conclusions}
In this work, a gradient-flow framework with a suitable time-stepping scheme and automatic tolerance adaptation for the solution of optimization problems with rather general constraints has been presented. Such a framework has then been applied to the problem of preliminary optimal design of a wave energy park.

From the point of view of numerical algorithms, the presented framework could be used as the local optimization step in a hybrid local-global method (see e.g. \cite{Wales_1997, Requena2006, dangelo2021} and the application to wave energy \cite{Neshat_2020}).  Moreover, we have considered the simplest adaptive Runge-Kutta scheme. Other promising choices are multistep methods, of which Nesterov's optimization algorithm can be considered a particular case \cite{Scieur2017}, and Runge-Kutta-Chebyshev methods \cite{Eftekhari_2020}, which introduce additional stages to improve stability and thus allow the use of very large time steps.

Concerning the application, the main limitation of this work is the restriction to cylindrical devices, which allows the analytic computation of derivatives. The presented framework is based on interaction theory and may be extended to more general shapes through the strategy introduced in \cite{McNatt2015}. More advanced control strategies could additionally be taken into account.


\section*{Acknowledgments}
Development of the \texttt{floatcyl} code commenced under the supervision of professor Giuseppe Passoni and contains contributions by Jacopo Gallizioli.
All authors are members of the GNCS Indam group. The present research is part of the activities of
“Dipartimento di Eccellenza 2023-2027."
\clearpage
\bibliographystyle{siamplain}
\bibliography{paper}

\end{document}